\documentclass{amsart}

\input xy
\xyoption{all}

\usepackage{geometry}
\usepackage{amsmath}
\usepackage{amssymb}
\usepackage{hyperref}
\usepackage{cite}
\usepackage[final]{showkeys}
\usepackage{hyperref}
\usepackage{setspace}
\usepackage{amsthm}  
\usepackage[pagewise,displaymath,mathlines]{lineno}

\newtheorem{same}{This should never appear}[section]
\newtheorem{defin}[same]{Definition}

\newtheorem{remark}[same]{Remark}
\newtheorem{theorem}[same]{Theorem}

\newtheorem{lemma}[same]{Lemma}
\newtheorem{fact}[same]{Fact}
\newtheorem{question}[same]{Question}
\newtheorem{cor}[same]{Corollary}
\newtheorem{prop}[same]{Proposition}

\newbox\noforkbox \newdimen\forklinewidth
\setbox0\hbox{$\textstyle\smile$}
\setbox1\hbox to \wd0{\hfil\vrule width \forklinewidth depth-2pt
 height 10pt \hfil}
\setbox\noforkbox\hbox{\lower 2pt\box1\lower 2pt\box0\relax}
\def\unionstick{\mathop{\copy\noforkbox}\limits}

\def\nonfork_#1{\unionstick_{\textstyle #1}}

\setbox0\hbox{$\textstyle\smile$}
\setbox1\hbox to \wd0{\hfil{\sl /\/}\hfil}
\setbox2\hbox to \wd0{\hfil\vrule height 10pt depth -2pt width
              \forklinewidth\hfil}
\newbox\doesforkbox
\setbox\doesforkbox\hbox{\lower 2pt\box1 \lower 2pt\box2\lower2pt\box0\relax}
\def\nunionstick{\mathop{\copy\doesforkbox}\limits}

\def\fork_#1{\nunionstick_{\textstyle #1}}

\newcommand{\ON}{\textrm{\textbf{ON}}}

\newcommand{\ba}{\bold{a}}
\newcommand{\bb}{\bold{b}}
\newcommand{\bc}{\bold{c}}

\newcommand{\bm}{\bold{m}}
\newcommand{\bn}{\bold{n}}

\newcommand{\bx}{\bold{x}}

\newcommand{\bz}{\bold{z}}
\newcommand{\dom}{\textrm{dom }}
\newcommand{\Th}{Th}

\newcommand{\rest}{\upharpoonright}
\newcommand{\id}{\textrm{id}}

\newcommand{\rank}{\text{rank }}

\newcommand{\crit}{\text{crit }}
\newcommand{\otp}{\text{otp}}

\newcommand{\bs}{\bold{s}}

 \newcommand{\footnotei}[1]{}

\newcommand{\seq}[1]{\langle #1 \rangle}

\newcommand{\bL}{\mathbb{L}}

\newcommand{\bX}{\bold{X}}

\newcommand{\cC}{\mathcal{C}}
\newcommand{\cF}{\mathcal{F}}
\newcommand{\cL}{\mathcal{L}}
\newcommand{\cM}{\mathcal{M}}
\newcommand{\cN}{\mathcal{N}}
\newcommand{\cP}{\mathcal{P}}

\newcommand{\comment}[1]{}

\title{Model theoretic characterizations of large cardinals}
\author{Will Boney}
\email{wboney@math.harvard.edu}
\address{Department of Mathematics, Harvard University, Cambridge, MA, USA}
\date{\today\\This material is based upon work done while
the author was supported by the National Science Foundation under Grant No. DMS-1402191.}

\begin{document}
\maketitle

\begin{abstract}
We consider compactness characterizations of large cardinals.  Based on results of Benda \cite{b-sccomp}, we study compactness for omitting types in various logics.  In $\bL_{\kappa, \kappa}$, this allows us to characterize any large cardinal defined in terms of normal ultrafilters, and we also analyze second-order and sort logic.  In particular, we give a compactness for omitting types characterization of huge cardinals, which have consistency strength beyond Vop\v{e}nka's Principle.
\end{abstract}

\section{Introduction}

Large cardinals typically have many equivalent formulations: elementary embeddings, ultrafilters or systems of ultrafilters, combinatorial properties, etc.  We investigate various characterizations in terms of logical compactness.  These formulations have a long history.  Weakly and strongly compact cardinals were first isolated with generalizations of the compactness theorem to infinitary languages.
\begin{fact}[\cite{t-compactness}]\
\begin{enumerate}
	\item $\kappa > \omega$ is weakly compact iff every $<\kappa$-satisfiable theory of size $\kappa$ in $\bL_{\kappa, \kappa}$ is satisfiable.
	\item $\kappa > \omega$ is strongly compact iff every $<\kappa$-satisfiable theory in $\bL_{\kappa, \kappa}$ is satisfiable.
\end{enumerate}
\end{fact}
Measurable cardinals also have such a characterization, this time in terms of \emph{chain} compactness of $\bL_{\kappa, \kappa}$.  This result is interesting because it seems to have been well-known in the past (evidenced by the fact that it appears as an exercise in Chang and Keisler's {\bf Model Theory} \cite[Exercise 4.2.6]{changkeisler}), but seems to have fallen out of common knowledge even among researchers working in the intersection of set theory and model theory (at least among the younger generation)\footnote{This is based on the author's personal impressions.  Although the statement seems forgotten, the proof is standard: if $T = \cup_{\alpha < \kappa} T_\alpha$ and $M_\alpha \vDash T_\alpha$, then set $M_\kappa : = \prod M_\alpha /U$ for any $\kappa$-complete, nonprincipal ultrafilter $U$ on $\kappa$.  \L o\'{s}' Theorem for $\bL_{\kappa, \kappa}$ implies $M_\kappa \vDash T$.}.

\begin{fact} \label{meas-fact}
$\kappa$ is measurable iff every theory $T \subset \bL_{\kappa, \kappa}$ that can be written as a union of an increasing $\kappa$-sequence of satisfiable theories is itself satisfiable.
\end{fact}

Magidor \cite[Theorem 4]{m-roleof} showed that extendible cardinals are the compactness cardinals of second-order logic, and Makowsky \cite{m-vopcomp} gives  an over-arching result that Vop\v{e}nka's Principle is equivalent to the existence of a compactness cardinal for every logic (see Fact \ref{stavi-fact} below).  This seems to situate Vop\v{e}nka's Principle as an upper bound to the strength of cardinals that can be reached by compactness characterizations.

However, this is not the case.  Instead, a new style of compactness is needed, which we call \emph{compactness for omitting types} (Definition \ref{cto-def}).  Recall that a type $p(x)$ in a logic $\cL$ is a collection of $\cL$-formulas in free variable $x$.  A model $M$ realizes $p$ if there is $a \in M$ realizing every formula in $p$, and omits $p$ if it does not realize it.  Explicitly, this means that for every $a \in M$, there is $\phi(x) \in p$ such that $M \vDash \neg \phi(a)$.  Although realizing a type in $\cL$ can be coded in the same logic by adding a new constant, omitting a type is much more difficult.  Thus, while in first-order logic types can be realized using a simple compactness argument, finding models omitting a type is much more difficult.  This is the inspiration for Gerald Sacks' remark\footnote{Several colleagues have suggested that Sacks is quoting himself here, as he considers himself a well-known recursion theorist and not a well-known model theorist.} \cite[p. 64]{s-saturated} ``A not well-known model theorist once remarked: `Any fool can realize a type, but it takes a model theorist to omit one.' ''

Compactness for omitting types was first (and seemingly uniquely) used by Benda to characterize supercompact cardinals in $\bL_{\kappa, \kappa}$.   Tragically, Benda's fantastic result is even less well-known than the characterization of measurable cardinals.  In the almost 40 years since its publication, the publisher reports zero citations of Benda's paper.

\begin{theorem}[{\cite[Theorem 1]{b-sccomp}}] \label{benda-fact}
Let $\kappa \leq \lambda$.  $\kappa$ is $\lambda$-supercompact iff for every $\bL_{\kappa, \kappa}$-theory $T$ and type $p(x, y) = \{\phi_i(x, y) \mid i < \lambda\}$, if there are club-many $s \in \cP_\kappa \lambda$ such that there is a model of 
$$T \cup \left\{ \exists x \left( \bigwedge_{i \in s} \exists y \phi_i(x,y) \wedge \neg \exists y \bigwedge_{i \in s} \phi_i(x, y)\right) \right\}$$
then there is a model of 
$$T \cup \left\{ \exists x \left( \bigwedge_{i < \lambda} \exists y \phi_i(x,y) \wedge \neg \exists y \bigwedge_{i < \lambda} \phi_i(x, y)\right) \right\}$$
\end{theorem}

The final model omits the type $q(y) = \{\phi_i(a, y) \mid i < \lambda\}$, where $a$ is the witness to the existential.  Phrased in these terms, this property says that if every small part of a type can be omitted, then the whole type can be omitted.  One complicating factor is that monotonicity for type omission works in the reverse direction as for theory satisfaction: larger types are \emph{easier} to omit since they contain more formulas.  This makes Benda's result somewhat awkward to phrase as he fixes the theory.

Our phrasing of compactness for omitting types varies from Benda's formulation in two key ways.  First, we also allow the theory to be broken into smaller pieces, which makes the phrasing more natural (at least from the author's perspective).  Second, and more crucially, we look at other index sets (or templates) than $\cP_\kappa \lambda$.  This allows us to capture many more large cardinals than just supercompacts.

Other model-theoretic properties have also been used to characterize large cardinals, mainly in the area of reflection properties and the existence of L\"{o}wenheim-Skolem-Tarski numbers.  Magidor \cite[Theorem 2]{m-roleof} characterizes supercompacts as the L\"{o}wenheim-Skolem-Tarski numbers for second-order, and Magidor and V\"{a}\"{a}n\"{a}nen \cite{mv-lst} explore the possibilities surrounding the L\"{o}wenheim-Skolem-Tarski numbers of various fragments of second-order logic.  Bagaria and V\"{a}\"{a}n\"{a}nen \cite{bv-symb} connect structural reflection properties and L\"{o}wenheim-Skolem-Tarski numbers through V\"{a}\"{a}n\"{a}nen's notion of symbiosis.  Of course, Chang's Conjecture has long been known to have large cardinal strength (see \cite[Section 7.3]{changkeisler}).

Section \ref{prelim-sec} fixes our notation and gives some basic results.  Section \ref{fo-sec} establishes the main definitions of compactness for type omission and applies it to the logics $\bL_{\kappa, \kappa}$.  The main result in this section is Theorem \ref{framework-thm}.  Section \ref{second-order-sec} examines type-omitting compactness for higher-order (Theorem \ref{second-thm}) and sort (Theorem \ref{sort-thm}) logics.  This section also deals with compactness characterizations of some other cardinals (e.g., strong) and discusses the notion of elementary substructure in second-order logic.  Section \ref{extot-sec} discusses extenders and type omission around the following question: we give characterizations of large cardinals with various logics, from $\bL_{\kappa,\kappa}$ to the all-powerful $\bL^{s, \Sigma_n}$.  However, $\bL_{\kappa, \kappa}$ is already able to work its way up the large cardinal hierarchy, including $n$-hugeness (Corollary \ref{main-cor}) and rank-into-rank (Section \ref{rir-sec}).  Are these more powerful logics necessary?  In other words, can we characterize all model-theoreticaly characaterizable large cardinals (extendible, etc.) by some property of $\bL_{\kappa, \kappa}$, or is the use of stronger logics necessary to pin down certain cardinals?  We focus on strong cardinals, and give a theory and collection of types such that the ability to find a model of the theory omitting the types is equivalent to $\kappa$ being $\lambda$ strong.  As close as this seems to an $\bL_{\kappa, \kappa}$ (or rather $\bL_{\kappa, \omega}(Q^{WF})$) characterization of strong cardinals, we still lack a general compactness for type omission for this case.

Preliminary results along these lines were first presented at the Workshop on Set-Theoretical Aspects of the Model Theory of Strong Logics hosted by the Centre de Recerca Matem\`{a}tica in 2016, and I'd like to thank many of the participants for helpful conversations, especially Jouko V\"{a}\"{a}n\"{a}nen for discussions about sort logic.  I'd also like to thank Gabriel Goldberg for helpful discussions regarding the strength of huge-for-$\bL^2_{\kappa, \kappa}$ cardinals, and Adrian Mathias and Sebastien Vasey for comments on a preliminary draft.  I would also like to thank the anonymous referee for helping to improve the paper.

\section{Preliminaries} \label{prelim-sec}

We begin with an informal introduction to the logics used.  The large cardinals notions are standard; consult Kanamori \cite{kanamori} or the locally given citation for detail.  We introduce some new large cardinal notions, typically naming them and defining them in the statement of a result: see Corollary \ref{sc-clear-cor}, Proposition \ref{huge-second-prop},  and Theorems \ref{wc-second-thm} and \ref{ns-strong-thm}.


$\bL_{\omega, \omega}$ is the standard, elementary first-order logic.

$\bL_{\lambda, \kappa}$ augments $\bL_{\omega, \omega}$ by allowing
\begin{itemize}
	\item conjunctions of $<\lambda$-many formulas that together contain $<\kappa$-many free variables;
	\item $<\kappa$-ary functions and relations in the language; and
	\item universal and existential quantification over $<\kappa$-many variables at once.
\end{itemize}
We typically restrict to $\lambda \geq \kappa$, both regular.

$\bL^2 = \bL^2_{\omega,\omega}$ is second-order logic, which extends $\bL_{\omega, \omega}$ by allowing quantification over subsets of finite cartesian powers of the universe and has an atomic `membership' relation.  The standard interpretation of the second-order quantifiers is quantification over \emph{all} subsets of (finite cartesian powers of) the universe, but an important concept is the nonstandard Henkin models $(M, P, E)$, where $M$ is a $\tau$-structure, $E \subset M\times P$ is an extensional relation, and $P$ represents a collection of subsets that the second-order quantifiers can range over.  The class of Henkin models of a second-order theory reduces to the models of a sorted $\bL_{\omega, \omega}$-theory, but we will still find use for this definition in Definition \ref{full-to-def} and the characterization of strong cardinals.  There is a second-order sentence $\Psi$ that says a Henkin model is standard:
$$\forall X \subset M \exists x \in P \forall y \in M \left( y \in X \leftrightarrow y E x\right)$$
We can also introduce higher-order variants $\bL^n$, but these are all codeable in $\bL^2$, e.g., a Henkin model for $\bL^3$ is $(M, P, Q, E)$, where $P$ represents the subsets of $M$ and $Q$ represents the subsets of $M$.  Then second-order can express the standardness of this model for third-order logic by using $\Psi$ and a copy of $\Psi$ for $P$ and $Q$.  The $\bL^n$ are coded similarly, although we must use $\bL^2_{|\alpha|^+,\omega}$ for $\bL^\alpha$ when $\alpha$ is infinite.  Thus, it suffices to talk about second order logic.

Additionally, when dealing with second-order logic, we allow the language to include functions and relations  whose domain and range include the second-order part of the model.  Given such a second-order language $\tau$, we describe it as consisting of a strictly first-order part (the functions and relations only using first-order inputs) and a second-order part (the functions and relations that use first- and second-order inputs).

We want to distinguish a very helpful second-order sentence $\Phi$ that will allow us to make a statement about the type of model we are in.  The construction of this statement is from Magidor \cite{m-roleof}, appearing in the proofs of his Theorems 2 and 4.  However, it is useful to extract the construction so we can make explicit reference to it.

\begin{fact}\label{magidor-phi-fact}
There is a second-order sentence $\Phi$ in the language with a single binary relation $E$ such that $(A, E) \vDash \Phi$ iff $(A, E) \cong (V_\alpha, \in)$ for some limit ordinal $\alpha$.  Moreover, $\Phi$ is a $\Pi^1_1$ sentence.
\end{fact}

{\bf Proof:} There is a first-order sentence $\phi(x,y)$ that is true iff $x$ is an ordinal and $y$ is a copy of $V_x$ (as computed in the model).  Then $\Phi$ is the conjunction of the statements (note that we use `$E$' for the relation in the language and `$\in$' for the logical notion of membership between a first- and second-order variable)
\begin{itemize}
	\item $E$ is well-founded: $\forall X \exists x \forall y \neg(y Ex \wedge y \in x)$;
	\item $E$ is extensional: $\forall x,y \left(x=y  \leftrightarrow \forall z (z E x \leftrightarrow z E y)\right)$;
	\item Every ordinal is in a $V_\alpha$: $\forall x \left(On(x) \to \exists y \phi(x,y)\right)$; and 
	\item Every subset of an element is represented in the model: $\forall x \forall X \left( \forall y (y \in X \to y Ex) \to \exists z \forall y (y \in X \leftrightarrow y E z)\right)$.
	
\end{itemize}
\hfill \dag\\

$\bL(Q^{WF})$ is $\bL_{\omega, \omega}$ augmented by the quantifier $Q^{WF}$ that takes in two free variables and so $Q^{WF} xy\phi(x, y, \bz)$ is true iff there is no infinite sequence $\{x_n \mid n < \omega\}$ such that $\phi(x_{n+1}, x_n, \bz)$ holds for all $n < \omega$; that is, $\phi(x, y, \bz)$ defines a well-founded relation.  Note that, in models of some choice, $Q^{WF}$ is both $\bL_{\omega_1, \omega_1}$ and $\bL^2$ expressible.  However, it will be useful to have it, e.g., in Theorem \ref{strong-char-thm}.

Finally, sort logic $\bL^s$ is a logic introduced by V\"{a}\"{a}n\"{a}nen \cite{v-sort-intro}.  This augments second-order logic by adding sort quantifiers $\exists^\sim, \forall^\sim$ where $\exists^\sim X \phi(X, \bx)$ is true in a structure $M$ iff there is a set $X$ (any set, not just a subset of the universe of $M$) such that $\phi(X, \bx)$ is true.  Sort logic is very powerful because it allows one to access a large range of information regardless of the language of the initial structure.  For instance, one can easily write down a formula $\Phi$ whose truth in any structure implies the existence of an inaccessible cardinal.  V\"{a}\"{a}n\"{a}nen discusses its use as a foundation of mathematics in \cite{v-sort-found}.  Since sort logic involves satisfaction of formulas in $V$, for definability of truth reasons, we must restrict to the logics $\bL^{s, \Sigma_n}$, where $\bL^{s, \Sigma_n}$ consists only of the formulas of sort logic that are of $\Sigma_n$ complexity when looking at the quantifiers over sorts.

Finally, all of these logics can be combined in the expected way, e.g., $\bL^2_{\kappa, \kappa}$.  We often take the union of two logics, e.g., $\bL^2 \cup \bL_{\kappa, \omega}$ is the logic whose formulas are in $\bL^2$ or $\bL_{\kappa, \omega}$; however, no second-order quantifier or variable can appear in any formula with an infinite conjunction, which separates it from $\bL^2_{\kappa, \omega}$.  We typically use boldface $\bL$ when discussing a particular logic and script $\cL$ when discussing an abstract logic.

For a logic $\cL$ and a language $\tau$, an $\cL(\tau)$-theory $T$ is a collection of sentences (formulas with no free variables) of $\cL(\tau)$.  An $\cL(\tau)$-type $p(x)$ in $x$ is a collection of formulas from $\cL(\tau)$ all of whose free variables are at most $x$.\footnote{Types in single variables suffice for the various characterizations in the paper, but they also extend to types of arity $<\kappa$.}  A type $p(x)$ is realized in a $\tau$-structure $M$ iff there is an element of the model that satisfies every formula in it and a type is \emph{omitted} precisely when it is not realized.  Note that the ``monotonicity of type omission" works the opposite way as theories: if $p(x) \subset q(x)$ are both types, then it is easier to omit $q$ than $p$.  We will often refer to filtrations of a theory $T$.  This means there is some ambient partial order $(I, \subset)$ and a collection of theories $\{T_s \mid s \in I\}$ such that $T = {\textstyle \bigcup}_{s\in I} T_s$ and $s \subset t$ implies $T_s \subset T_t$.

In general, we are agnostic about how one codes these logics as sets, except to insist that it is done in a reasonable way, e.g., $\tau$ is coded as a set of rank $|\tau|+\omega$, $\bL_{\kappa, \kappa}(\tau) \subset V_{\kappa +|\tau|}$, etc.  This gives us two nice facts about the interaction between languages $\tau$ and elementary embeddings $j:V \to \cM$ (or $V_\alpha \to V_\beta$, etc.) with $\crit j =\kappa$:
\begin{itemize}
	\item if $\tau$ is made up of $< \kappa$-ary functions and relations, then $j"\tau$ and $\tau$ are just renamings of each other; and
	\item if $\phi \in \bL^{s, \Sigma_n}_{\kappa, \kappa}(\tau)$, then $j(\phi) \in \bL^{s, \Sigma_n}_{\kappa, \kappa}(j"\tau)$.
\end{itemize}
This means that, when searching for a model of $T$, it will suffice to find a model of $j"T$, which is a theory in the same logic and an isomorphic language.

Given an inner model\footnote{In an unfortunate collision of notation, $M$ is commonly used for both inner models in set theory and for $\tau$-structures in model theory.  Owing to my model-theoretic roots, this paper uses standard $M$ for $\tau$-structures and script $\cM$ for inner models.} $\cM$ (or some $V_\alpha$), we collect some facts about when $\cM$ is correct about various logics.  That is, the statement ``$M$ is a $\tau$-structure" is absolute from $\cM$ to $V$ and we want to know when the same holds of ``$\phi$ is an $\cL(\tau)$-formula and $M\vDash_{\cL} \phi$."
\begin{itemize}
	\item $\cM$ is correct about the logic $\bL_{ON^{\cM}, \omega}(Q^{WF})$.
	\item If ${}^{<\kappa}\cM \subset \cM$, then $\cM$ is correct about the logic $\bL_{\kappa, \kappa}$.
	\item If $\cP(A) \in \cM$, then $\cM$ is correct about $\bL^2$ for structures with universe $A$.
	\item If $\cM \prec_{\Sigma_n} V$, then $\cM$ is correct about $\bL^{s, \Sigma_n}$.
\end{itemize}

As a warm-up, note that any compactness involving an extension of $\bL(Q^{WF})$ will entail the existence of large cardinals.  Fixing $\kappa$, if 
$$T:=ED_{\bL_{\omega, \omega}}\left(V_{\kappa+1}, \in, x\right)_{x \in V_{\kappa+1}} \cup \{c_\alpha < c < c_\kappa\mid \alpha < \kappa\}\cup\{Q^{WF} xy (x \in y)\}$$ 
is satisfiable, then there is a non-surjective elementary embedding $j:V_{\kappa+1} \to \cM$ to a well-founded structure with $\crit j \leq \kappa$.  Standard results imply that $\crit j$ must be measurable.  Moreover, $T$ is `locally satisfiable' in the sense that, if $T_0 \subset T$ does not contain constants for elements with ranks unbounded in $\kappa$, then $V_{\kappa+1}$ can be made a model of $T_0$ by adding constants.

\section{Type-omitting compactness in $\bL_{\infty, \infty}$} \label{fo-sec}

We introduce some basic definitions that will be used in each of our characterizations.  The notion of containing a strong $\kappa$-club becomes very important, and we discuss that concept after the definition.

\begin{defin} \label{basic-def}
Let $\kappa \leq \lambda$ and $I \subset \cP(\lambda)$.
\begin{enumerate}
	\item $I$ is \emph{$\kappa$-robust} iff for every $\alpha < \lambda$, $[\alpha]_I:=\{ s \in I \mid \alpha \in s\} \neq\emptyset$ and $I\subset\{s : |s \cap \kappa|<\kappa\}$.
	\item $C \subset I$ \emph{contains a strong $\kappa$-club} iff there is a function $F:[\lambda]^2 \to \cP_\kappa \lambda$ such that
	$$C(F):= \{s \in I\mid s \text{ is infinite and, for all }x, y, \in s, F(x, y) \subset s\} \subset C$$
\end{enumerate}
Let $U$ be an ultrafilter on $I$.
\begin{enumerate}
	\item[(3)] $U$ is \emph{$\mu$-complete} iff for all $\alpha < \mu$ and $\{X_\beta \in U \mid \beta < \alpha\}$, we have $\bigcap_{\beta<\alpha} X_\beta \in U$.
	\item[(4)] $U$ is \emph{fine} iff for all $\alpha \in \lambda$, $[\alpha]_I \in U$.
	\item[(5)] $U$ is \emph{normal} iff for all $F:I \to \lambda$ such that $\{s\in I\mid F(s) \in s\} \in U$, there is $\alpha_0<\lambda$ such that $\{s \in I \mid F(s) = \alpha_0\} \in U$. 
\end{enumerate}
\end{defin}

The conditions of $\kappa$-robustness are intended to make sure that $I$ includes enough sets so that the notion of a ``$\kappa$-complete, normal, fine ultrafilter on $I$" makes sense and is possible.  In particular, any such ultrafilter $U$ will be characterized by an elementary embedding $j_U$ with $\crit j_U = \kappa$; this implies that $\{s \in I : |s \cap \kappa |<\kappa\} \in U$.

We define the notion of `contains a strong $\kappa$-club' without defining the notion of a strong $\kappa$-club.  However, the  first two items of Fact \ref{right-gen-fact} show that this notion correctly generalizes the notion of containing a club from $\kappa$ and $\cP_\kappa \lambda$.  Moreover, the third item shows that a generalization that replaces $\cP_\kappa\lambda$ with a different set does not work.

\begin{fact}\label{right-gen-fact} \
\begin{enumerate}
	\item If $I = \cP_\kappa\lambda$, then containing a strong $\kappa$-club is equivalent to containing a club.
	\item Fix a $\kappa$-robust $I \subset \cP(\lambda)$.  If $U$ is a $\kappa$-complete, fine, normal ultrafilter on $I$, then it extends the contains a strong $\kappa$-club filter, that is, $C(F) \in U$ for all $F:[\lambda]^2\to \cP_\kappa \lambda$. 
	\item If $I = [\lambda]^{\kappa}$ and $U$ is a $\kappa$-complete, fine, normal ultrafilter on $I$, then there is no $s \in [\lambda]^\kappa$ such that
	$$[s]_I:=\{t \in I \mid s \subset t\} \in U$$
\end{enumerate}
\end{fact}

{\bf Proof:} (1) is a result of Menas, see \cite[Proposition 25.3]{kanamori}.

For (2), fix $F: [\lambda]^2 \to \cP_\kappa \lambda$ and suppose that $C(F) \not \in U$.  For each $s \in I- C(F)$, there are $\alpha_s < \beta_s \in s$ such that $F(\alpha_s, \beta_s) \not\subset s$.  By applying normality twice, there is some $Z \in U$ and $\alpha_* < \beta_*$ such that, for all $s \in Z$, $\alpha_s = \alpha_*$ and $\beta_s = \beta_*$.  By the $\kappa$-completeness and fineness of $U$, we have that $[f(\alpha_*, \beta_*)]_{I} \in U$.  Thus, there is $t \in Z \cap [F(\alpha_*, \beta_*)]_{I}$; however, this is a contradiction.

For (3), given such a $U$, build the elementary embedding $j_U: V \to \cM \cong \prod V/U$.  Let $s \in [\lambda]^\kappa$.  Then, for $X \subset [\lambda]^\kappa$, we have that $X \in U$ iff $j_U''\lambda \in j_U(X)$.  However, since $|s|=\kappa = \crit j_U$, there is some $\alpha \in j_U(s) - j_U"\lambda$.  In particular, this means that
$$j_U''\lambda \not\in \{ t \subset j(I) \mid j_U(s) \subset t\} = j_U\left([s]_I\right)$$
\hfill\dag\\

We are interested in model-theoretic conditions that guarantee the existence of a fine, normal, $\kappa$-complete ultrafilter on some $\kappa$-robust $I$.  Recall that, from Kunen's proof of the inconsistency of Reinhardt cardinals, every countably complete ultrafilter on $\cP(\lambda)$ must concentrate on $\cP_{\mu} \lambda$ for some $\mu \leq \lambda$.  If $\mu$ is strictly larger than the completeness of the ultrafilter, we will need the following technical condition.  In practice, the set $X$ will always be a theory.

\begin{defin}
Let $I \subset \cP(\lambda)$ and $X$ be a set that is filtrated as an increasing union of $\{X_s \mid s \in I\}$.  Then we say this filtration \emph{respects the index} iff there is a collection $\{X^\alpha \mid \alpha \in \lambda \}$ such that, for each $s \in I$, $X_s = \bigcup_{\alpha \in s} X^\alpha$.
\end{defin}

This condition says that the filtration at $s \in I$ is just determined by the elements of $s$.  Note this condition is trivially satisfied when $I \subset \cP_\kappa \lambda$, but will be important in certain cases, e.g., to characterize huge cardinals (see Corollary \ref{main-cor}.(\ref{h})).  There, it will guarantee that if $\phi \in \cup_{s \in [\lambda]^\kappa} T_s$, then there are a large number of $s \in I$ such that $\phi \in T_s$.

The main concept of this section is the following:

\begin{defin}\label{cto-def}
Let $\cL$ be a logic (in the sense of Barwise \cite{b-abstract} or taken without formal definition), $\kappa \leq \lambda$, and $I \subset \cP(\lambda)$ be $\kappa$-robust.  Then we say that \emph{$\cL$ is $I$-$\kappa$-compact for type omission} iff the following holds:

Suppose that we have
\begin{itemize}
	\item a language $\tau$;
	\item a $\cL(\tau)$-theory $T$ that can be written as an increasing union of $\{T_s :s \in I\}$ that respects the index; and
	\item a collection of types $\{p^a(x) \mid a \in A\}$ of size $\lambda$ (for an arbitrary set $A$), where each type comes with an enumeration $p^a(x) = \{\phi^a_\alpha(x) : \alpha < \lambda\}$ and, for $s \subset \lambda$, we set $p^a_s(x):=\{\phi^a_\alpha(x):\alpha \in s\}$.
\end{itemize}
If
\begin{equation} \label{cto-set} \tag{CTO}
\left\{ s \in I \mid T_s \text{ has a model omitting each type in }\{p^a_s(x) \mid a \in A\}\right\}
\end{equation}
contains a strong $\kappa$-club, then there is a model of $T$ omitting each type in $\{ p^a(x) \mid a \in A\}$.
\end{defin}

This definition is technical, so we will try to unpack it.  A compactness for type omissions result should say that if we want to find a model of a theory $T$ omitting a type $p$, it should suffice to find a model of all of the small fragments of $T$ omitting the small fragments of $p$.  When there is no $p$, the monotonicity of satisfying theories makes looser statements of this possible.  However, the monotonicity of satisfying theories and omitting types work in opposite directions (smaller theories are easier to satisfy, while smaller types are harder to omit), so much more care has to be taken.

The $s \in I$ index allows us to associate particular fragments of the theory with particular fragments of the type.  In this sense, wether or not  the set (\ref{cto-set}) contains a strong $\kappa$-club depends very strongly on the particular filtration of $T$ and enumerations of the $p^a$ that are chosen.  Additionally, because the monotonicities work in opposite directions, given $s \subset t \in I$, wether or not $s$ and $t$ are in (\ref{cto-set}) is independent.  As explained in Remark \ref{ultrafilter-rem}, we want the set (\ref{cto-set}) to be in whatever fine, normal, $\kappa$-complete ultrafilter our assumptions give us, and the results of Fact \ref{right-gen-fact} suggest that `contains a strong $\kappa$-club' is the right stand in for this notion.

The following gives a framework result for linking compactness for type omission for $\bL_{\kappa, \lambda}$ to various large cardinal notions.

\begin{theorem} \label{framework-thm}
Let $\kappa \leq \lambda$, and $I \subset \cP(\lambda)$ be $\kappa$-robust.  The following are equivalent:
\begin{enumerate}
	\item $\bL_{\kappa, \omega}$ is $I$-$\kappa$-compact for type omission.
	\item $\bL_{\kappa, \kappa}$ is $I$-$\kappa$-compact for type omission.
	\item \label{ult} There is a fine, normal, maximally $\kappa$-complete ultrafilter on $\cP(\lambda)$ concentrating on $I$; `maximally $\kappa$-complete' means $\kappa$-complete and not $\kappa^+$-complete.
	\item \label{emb} There is an elementary $j:V \to \cM$ with $\crit j = \kappa$ and $j"\lambda \in \cM \cap j(I)$.
\end{enumerate}
Moreover, the first $\mu$ such that $\bL_{\mu, \omega}$ is $I$-$\kappa$-compact for type omission is the first $\mu$ with a fine, normal, $\mu$-complete ultrafilter on $I$.
\end{theorem}

{\bf Proof:} The equivalence of $(3)$ and $(4)$ is straightforward from standard methods, and $(2)$ implies $(1)$ is obvious because $\bL_{\kappa, \omega} \subset \bL_{\kappa, \kappa}$.\\

$(4) \to (2)$: Suppose we have a set-up for $I$-$\kappa$-compact type omission: a theory $T$ in $\bL_{\kappa, \kappa}(\tau)$ filtrated as $T_s$ for $s \in I$ and a collection $\{p^a:a \in A\}$ of $\bL_{\kappa, \kappa}(\tau)$-types as in Definition \ref{cto-def}.  Set $\bar{T}$ and $\bar{p^a}$ to be the functions that take $s \in I$ to $T_s$ and $p_s^a$, respectively.  Similarly, set $\bar{M}$ to be the partial function that takes $s$ to $M_s$ that models $T_s$ and omits each $p^a_s(x)$ (when defined).

Since the set (\ref{cto-set}) contains a strong $\kappa$-club, there is a function $F:[\lambda]^2 \to \cP_\kappa\lambda$ that witnesses this.  Then $j(\bar{M})$ is a function with domain containing $j\left(C(F)\right)$.  We know that $j$ comes from an ultrafilter as in (3) and Fact \ref{right-gen-fact}.(2) implies that $C(F) \in U$.  Thus, $j"\lambda \in j\left(C(F)\right)$.

Then $\cM$ thinks that $j(\bar M)(j"\lambda)$ is a $j(\tau)$-structure that models $j(\bar T)(j"\lambda)$ that omits each $j(\bar p^a)(j"\lambda)$.  $\cM$ is correct about this, so $j(\bar M)(j"\lambda)$ models $j(\bar T)(j"\lambda)$ and omits each $j(\bar p^a)(j"\lambda)$ (in $V$). As discussed in Section \ref{prelim-sec}, $\tau$ and $j"\tau$ are essentially the same, in the sense that there is a canonical bijection between them--taking $F \in \tau$ to $j(F)\in \tau$--that respects arity; such a function is called a renaming.

{\bf Claim:} $j"T \subset j(\bar{T})(j"\lambda)$ and $j"p^a = j(\bar{p}^a)(j"\lambda)$ for each $a \in A$.\\
First, let $j(\phi) \in j"T$ and let $\{T_0^\beta:\beta \in \lambda\}$ witness that the filtration respects the index; this is the crucial use of this property.  Then we have that $\phi \in T_0^\alpha$ for some $\alpha \in \lambda$.  Then 
$$j(\phi) \in j(\bar{T_0})(\alpha) \subset j(\bar{T})(j"\lambda)$$
The proof that $j"p^a \subset j(\bar{p}^a)(j"\lambda)$ is similar.  For the other direction, let $\psi \in j(\bar{p}^a)(j"\lambda)$.  By elementarity, we have that $j(p^a) = \{\psi^a_\alpha(x): \alpha < j(\lambda)\}$ and $j(\bar{p}^a)(s) = \{\psi^a_\alpha(x):\alpha \in s\}$, where $\psi^a_{j(\beta)} = j(\phi^a_\beta)$.  This means that $\psi$ is of the form $j(\phi)$ for $\phi \in p^a$, which exactly says $\psi \in j"p^a$.  Here, we have exactly used the condition that $j"\lambda$ is in the model, and simply using a set that contained this would not work.

With the claim proved, we have produced a model $j(\bar{M})(j"\lambda)$ of $j"T$ that omits each $j"p^a$ for $a \in A$.  After applying the inverse of the canonical naming above, we have the desired model of $T$ that omits each $p^a$.\\

$(1) \to (3)$: We want to write down a theory $T$ and collection of types $\{p_F:F:I \to \lambda\}$ such that a model of $T$ omitting these types will code the desired ultrafilter.  Set $\tau = \{P, Q, E, c_X, d\}_{X \subset I}$ with $P$ and $Q$ unary predicates, $E \subset P \times Q$ a binary relation, and $c_X, d$ constants.  We look at the standard structure $M = \seq{I, \cP(I), \in, X}_{X \subset I}$ that has no interpretation for $d$.  Classic results tell us that finding an extension $N$ of $M$ that a new element $n$ (which will be the interpretation of $d$) codes an ultrafilter $U$ on $I$ by
$$X \in U \iff N \vDash n E c_X$$
The completeness of this ultrafilter comes from how strong the elementarity of the extension is, so we need to add conditions that make the fineness and normality hold.

Set $T_0 = \Th_{\bL_{\kappa, \omega}}(M)$ (although much less is necessary).  Set $T^\alpha := T_0 \cup \{``d \in c_{[\alpha]^I}"\}$ for $\alpha < \lambda$; $T_s := \cup_{\alpha \in s} T^\alpha$ for $s \in I$; and $T = \cup_{s \in I} T_s$.  For each function $F:I \to \lambda$, define 
\begin{eqnarray*}
X_F&:=&\{s \in I\mid F(s) \in s\}\\
X_{F, \alpha} &:=& \{s \in I \mid F(s) = \alpha\}\\
p_F(x) &:=& \left\{ x = d \wedge x E c_{X_F} \wedge \neg \left(x E c_{X_{F, \alpha}}\right) \mid \alpha < \lambda\right\}\\
\Gamma &:=& \{ p_F \mid F:I \to \lambda\}
\end{eqnarray*}
These types are built so that omitting the type $p_F$ means that the ultrafilter coming from $d$ will be normal with respect to the function $F$: if $F$ is regressive on a large set, then ``$dEc_{X_F}$'' will hold, so omitting $p_F$ requires that there is an $\alpha < \lambda$ such that ``$d E c_{X_{F, \alpha}}$'' holds, and $\alpha$ will then give the constant value of $F$ on a large set.

We will use compactness for type omission to find a model of $T$ omitting $\Gamma$.

{\bf Claim 1:} If there is a model of $T$ omitting $\Gamma$, then there is a fine, normal, $\kappa$-complete ultrafilter on $I$.\\
Let $N$ be this model.  Define $U$ on $I$ by 
$$X \in U \iff N\vDash d E c_X$$
It is straightforward to check that $U$ is a $\kappa$-complete ultrafilter on $I$.  For instance, given $\seq{X_\alpha \in U \mid \alpha < \mu < \kappa}$, we know that the following sentence is in $T$:
$$\forall x \left( \bigwedge_{\alpha < \mu} x E c_{X_\alpha} \to x E c_{\cap_{\alpha<\mu} X_\alpha}\right)$$
Thus, $N \vDash d E c_{\cap_{\alpha<\mu} X_\alpha}$.

Given $\alpha < \lambda$, by $\kappa$-robustness, there is some $s \in I$ such that $\alpha \in s$.  So $T$ entails that $d E c_{[\alpha]^I}$. This means that $U$ is fine.

For normality, if $F:I \to \lambda$ is regressive on a $U$-large set, then $N \vDash d E c_{X_F}$.  Since $N$ omits $p_F$, there is $\alpha < \lambda$ such that $X_{F, \alpha} \in U$, so $U$ is normal.

{\bf Claim 2:} For each $s \in I$, there is a model of $T_s$ omitting $\{p_{F, s} \mid F:I \to \lambda\}$.\\
Expand $M$ to $M_s$ by interpreting $d^{M_s} = s$.  This models $T_s$ since $s \in [\alpha]^I$ for each $\alpha \in s$ by definition.  Moreover, if there is $x \in M_s$ such that
$$M_s \vDash ``x = d \wedge x E X_F"$$
for some $F:I \to \lambda$, then $x = s$ and $F(s) \in s$, so there is $\alpha \in s$ such that $F(s) = \alpha$.  Thus
$$M_s \vDash x E X_{F, \alpha}$$
So $M_s$ omits each $p_{F, s}$.

By the $I$-$\kappa$-compactness for type omission, we are done.\\

The proof of the ``moreover" follows similarly.\hfill \dag\\

\begin{remark}\label{ultrafilter-rem}
We chose to show $(4)\to(2)$ above because it will more easily generalize to large cardinals characterized by the existence of extenders rather than ultrafilters (this is done in Section \ref{second-order-sec}).  However, the direct proof of $(3)\to(2)$ helps to emphasize the connection between normality of an ultrafilter and type omission, so we outline it here.  Suppose that $U$ is an ultrafilter as in Theorem \ref{framework-thm}.(3) and we have a language $\tau$, theory $T$, and set of types $\{p^a(x) \mid a \in A\}$ as in the hypothesis of Definition \ref{cto-def}, and let $X$ be the set in \ref{cto-set}.  For $s \in X$, there is a witnessing model $M_s$ and we can form the ultraproduct\footnote{Here we use a more liberal formalism for ultraproducts that allow certain structures to be empty by allowing the choice functions $f:I \to \bigcup M_s$ to be partial with $U$-large domain.} $\prod M_s/U$.  From \L o\'{s}' Theorem, the $\kappa$-completeness and fineness of $U$ means that $\prod M_s/U$ models $T$.

To prove type omission, suppose that $[f]_U \in \prod M_s/U$ and let $a \in A$.  Since $M_s$ omits $p^a_s(x)$, there is some $\alpha_s \in s$ such that
$$M_s \vDash \neg \phi^a_{\alpha_s}\left(f(s)\right)$$
The function $s \mapsto \alpha_s$ is regressive on the $U$-large set $\dom f$, so normality implies that there is a single $\alpha_*$ such that 
$$M_s \vDash \neg \phi^a_{\alpha_*}\left(f(s)\right)$$
for a $U$-large set of $s$.  Then, by \L o\'{s}' Theorem, 
$$\prod M_s/U \vDash \neg \phi^a_{\alpha_*}\left([f]_U\right)$$
so $[f]_U$ doesn't realize $p^a$.  Since $[f]_U$ and $a$ were arbitrary, the ultraproduct omits the types $\{p^a(x) \mid a \in A\}$.

\end{remark}

As an example of the moreover, $\bL_{\omega, \omega}$ satisfies Benda's supercompactness theorem iff $\cP_\omega \lambda$ carries a fine, normal measure that need not even be countably complete. 
Note that $\bL_{\omega, \omega}$ can never be $\cP_\omega \lambda$-$\omega$-compact for type omission: the $\max$ function shows that no fine ultrafilter on $\cP_\omega \lambda$ can be normal.

This general framework directly gives model theoretic characterizations of large cardinals that are characterized by normal ultrafilters.

\begin{cor} \label{main-cor}
For each numbered item below, all of its subitems are equivalent:
\begin{enumerate}
\item \label{m}\begin{enumerate}
	\item $\kappa$ is measurable.
	\item $\bL_{\kappa,\kappa}$ is $\cP_\kappa \kappa$-$\kappa$-compact for type omission.
\end{enumerate}
\item \label{sc} \begin{enumerate}
	\item $\kappa$ is $\lambda$-supercompact.
	\item $\bL_{\kappa,\kappa}$ is $\cP_\kappa \lambda$-$\kappa$-compact for type omission.
\end{enumerate}
\item \label{h} \begin{enumerate}
	\item $\kappa$ is huge at $\lambda$; that is, $\kappa$ is huge and there is $j:V\to \cM$ witnessing this with $j(\kappa)=\lambda$.
	\item $\bL_{\kappa,\kappa}$ is $[\lambda]^{\kappa}$-$\kappa$-compact for type omission.
\end{enumerate}
\item \begin{enumerate}
	\item $\kappa$ is $n$-huge at $\lambda_1, \dots, \lambda_n$.
	\item $\bL_{\kappa,\kappa}$ is $\{ s \subset \lambda : \forall i < n, |s \cap \lambda_{i+1}| = \lambda_i\}$-$\kappa$-compact for type omission.
\end{enumerate}
\end{enumerate}
\end{cor}

{\bf Proof:} The proof follows the standard characterizations of these notions in terms of normal ultrafilters (see \cite{kanamori}) and from Theorem \ref{framework-thm}.

For instance, in (\ref{h}), if $\kappa$ is huge at $\lambda$, then there $j:V \to \cM$ with $\crit j = \kappa$ with ${}^{j(\kappa)}\cM \subset \cM$ and $\lambda = j(\kappa)$.  In particular, this means that $j"\lambda \in \cM$ and belongs to $j\left([\lambda]^\kappa\right) = \left([j(\lambda)]^{j(\kappa)}\right)^{\cM}$.  Using Theorem \ref{framework-thm}.(4)$\to$(2), $\bL_{\kappa, \kappa}$ is $[\lambda]^\kappa$-$\kappa$-compact for type omission.

In the other direction, assume that $\bL_{\kappa,\kappa}$ is $[\lambda]^\kappa$-$\kappa$-compact for type omission.  By Theorem \ref{framework-thm}.(2)$\to$(3), there is a fine, normal $\kappa$-complete ultrafilter on $I$ concentrating on $[\lambda]^\kappa$.  By \cite[Theorem 24.8]{kanamori}, $\kappa$ is huge at $\lambda$.\hfill\dag\\

Item (\ref{sc}) is Benda's supercompactness theorem (Theorem \ref{benda-fact}). Item (\ref{m}) can be reformulated along the lines of Fact \ref{meas-fact}:
\begin{center}
If $T = \cup_{\alpha < \kappa} T_\alpha$ is an $\bL_{\kappa, \kappa}(\tau)$-theory and $p(x) = \{\phi_i(x) \mid i < \kappa\}$ is a type such that for every $\alpha < \kappa$, there is a model of $T_\alpha$ omitting $\{\phi_i (x) \mid i < \alpha\}$, then there is a model of $T$ omitting $p$.
\end{center}

This helps highlight the impact that the normality of an ultrafilter has on the resulting ultraproduct: if $U$ is a normal ultrafilter on $I \subset \cP(\lambda)$ and $\{M_s \mid s \in I\}$ are $\tau$-structures, then $\prod M_s/U$ omits any type $p = \{\phi_\alpha(x) \mid \alpha < \lambda\}$ such that $\{s \in I \mid M_s \text{ omits } p_s\} \in U$.

As mentioned above, any ultrafilter on $\cP(\lambda)$ concentrates on some $\cP_\mu\lambda$.  We can characterize when an ultrafilter exists on some $\cP_\mu\lambda$ with the following large cardinal notion:

\begin{cor}\label{sc-clear-cor}
Fix $\kappa \leq \mu \leq \lambda$.  The following are equivalent:
\begin{enumerate}
	\item $\kappa$ is $\lambda$-supercompact with $\mu$ clearing: there is $j:V \to \cM$ with $j"\lambda \in \cM$ and $j(\mu) > \lambda$.
	\item $\bL_{\kappa,\kappa}$ is $\cP_\mu \lambda$-$\kappa$-compact for type omission.
\end{enumerate}
\end{cor}

{\bf Proof:} This follows from Theorem \ref{framework-thm}: $j(\mu) > \lambda$ ensures that $j"\lambda \in j(\cP_\mu\lambda)$.\hfill\dag\\

This equivalence even holds for $\kappa=\omega$, where both conditions fail for all $\mu \leq \lambda$.

One feature of the compactness schema are that the theories are not required to have a specific size, but rather should be filtrated by a particular index set.  Note this is also true for strongly compact cardinals; that is, rather than characterizing $\lambda$-strongly compact cardinals as compactness cardinals for $\lambda$-sized theories in $\bL_{\kappa, \kappa}$, we can give the following.

\begin{prop}
$\kappa$ is $\lambda$-strongly compact iff any $\bL_{\kappa, \kappa}$-theory $T$ that can be filtrated as an increasing union of satisfiable theories indexed by $\cP_\kappa\lambda$ is itself satisfiable.
\end{prop}

This can even be extended to theories of proper class size.  Each item of Corollary \ref{main-cor} remains true if compactness for type omission is generalized to allow for the $T$ and the $T_s$ in Definition \ref{cto-def} to be definable proper classes and to allow for satisfaction by definable proper class structures.  The proof of Theorem \ref{framework-thm} goes through with this generalization.

\begin{remark}
$ED_\cL(V, \in, x)_{x \in V}$ is used below to denote the $\cL$-elementary diagram of the structure with universe $V$, a single binary relation $\in$, and a constant for each $x$ that is interpreted as $x$.  Formally, from Tarski's undefinability of truth, this is not a definable class when $\cL$ extends $\bL_{\omega, \omega}$.  Similarly, the statement that $j:V \to \cM$ is elementary is not definable.

However, e.g., \cite[Proposition 5.1.(c)]{kanamori} shows that, for embeddings between inner models, $\Sigma_1$-elementarity implies $\Sigma_n$-elementarity for every $n < \omega$.  Thus, mentions of elementary embeddings with domain $V$ can be replaced by $\Sigma_1$-elementarity.  Similarly, the full elementary diagram of $V$ could be replaced by its $\Sigma_1$-counterpart.  However, following set-theoretic convention, we continue to refer to the full elementary diagram.
\end{remark}

Armed with a class version of omitting types compactness, we can show equivalences directly between the model-theoretic characterizations and the elementary embedding characterizations without working through an ultrafilter characterization.  At each stage, we find a model $\cN$ of $ED_{\bL_{\kappa, \kappa}}(V, \in, x)_{x \in V}$ along with the sentences $\{c_\alpha < c < c_\kappa \mid \alpha < \kappa\}$, where $c$ is a new constant, and some other sentences.  Such a model is well-founded because it models the $\bL_{\omega_1, \omega_1}$-sentence that truthfully asserts well-foundedness, so we can take the Mostowski collapse $\pi: \cN \cong \cM$ with $\cM$ transitive.  Then $x \mapsto \pi(c_x)$ is an $\bL_{\kappa, \kappa}$-elementary embedding that necessarily sends $\alpha <\kappa$ to itself.  Moreover, the interpretation of $c$  guarantees that the critical point is at most $\kappa$ and the use of $\bL_{\kappa, \omega}$ guarantees the critical point is at least $\kappa$.  This is enough to show $\kappa$ is measurable and extra sentences to be satisfied and types to be omitted can be added to characterize the above large cardinal notions.

\begin{theorem}\label{elem-comp-thm}
For each numbered item below, all of its subitems are equivalent:
\begin{enumerate}
\item \begin{enumerate}
	\item $\kappa$ is measurable.
	\item There is a definable class model of
	$$ED_{\bL_{\kappa, \kappa}}(V, \in, x)_{x \in V} \cup \{c_\alpha < c < c_\kappa \mid \alpha < \kappa\}$$
\end{enumerate}
\item \begin{enumerate}
	\item $\kappa$ is $\lambda$-strongly compact.
	\item There is a model of
	$$ED_{\bL_{\kappa, \kappa}}(V, \in, x)_{x \in V} \cup \{c_\alpha < c < c_\kappa \mid \alpha < \kappa\} \cup \{c_\alpha \in d \wedge |d| < c_\kappa \mid \alpha < \lambda\}$$
\end{enumerate}
\item \begin{enumerate}
	\item $\kappa$ is $\lambda$-supercompact.
	\item There is a definable class model of
	$$ED_{\bL_{\kappa, \kappa}}(V, \in, x)_{x \in V} \cup \{c_\alpha < c < c_\kappa \mid \alpha < \kappa\} \cup \{c_\alpha \in d \wedge |d| < c_\kappa \mid \alpha < \lambda\}$$
	that omits
	$$p(x) = \{ x E d \wedge x \neq c_\alpha \mid \alpha < \lambda\}$$
\end{enumerate}
\item \begin{enumerate}
	\item $\kappa$ is $n$-huge at $\lambda_1, \dots, \lambda_n$.
	\item There is a definable class model of
	$$ED_{\bL_{\kappa, \kappa}}(V, \in, x)_{x \in V} \cup \{c_\alpha < c < c_\kappa \mid \alpha < \kappa\} \cup \{c_\alpha \in d_{i+1} \wedge |d_{i+1}| = c_{\lambda_i} \mid \alpha < \lambda_{i+1}, i < n\}$$
	that omits, for $i < n$,
	$$p_i(x) = \{ x E d_{i+1} \wedge x \neq c_\alpha \mid \alpha < \lambda_{i+1}\}$$
\end{enumerate}
\end{enumerate}
\end{theorem}

{\bf Proof:} We prove item (3); the rest of the proofs are similar.

First, suppose that $\kappa$ is $\lambda$-supercompact.  Then there is $j:V \to \cM$ with $\crit j=\kappa$, $j(\kappa)>\lambda$, and $j"\lambda \in \kappa$.  We claim that (the definable) $\cM$ is our model after expanding the vocabulary:
\begin{itemize}
	\item $c_x^\cM = j(x)$ for $x \in V$;
	\item $c^\cM = \kappa$; and
	\item $d^\cM = j"\lambda$.
\end{itemize}
The elementarity of $j$ and the closure under $\kappa$-sequences of $\cM$ imply that 
$$(V, \in, x)_{x \in V} \equiv_{\bL_{\kappa, \kappa}} \left(\prod V/U, E, [s \mapsto x]_U\right)_{x \in V} \cong \left(\cM, \in j(x)\right)_{x \in V}$$
The other sentences of the theory are true because $j(\alpha) < \kappa < j(\kappa)$ for $\alpha < \kappa$ precisely means $\crit j =\kappa$, and we know that $j(\alpha) \in j"\lambda$ for all $\alpha < \lambda$, with $|j"\lambda|=\lambda <j(\kappa)$.

Second, suppose that we have such a definable class model $(\cM, E, a_x, a, b)_{x \in V}$.  Well-foundedness is $\bL_{\omega_1, \omega}$-expressible\footnote{Note that the first-order Axiom of Foundation can hold in ill-founded classes.}:
$$\forall \seq{x_n:n<\omega} \neg \bigwedge_{n<\omega} x_{n+1} E x_n$$
Thus, $\cM$ is well-founded with respect to $E$.  By taking the Mostowski collapse, we can assume that $\cM$ is transitive and $E$ is $\in$.  Then define $j:V \to \cM$ by $j(x) = a_x$.  

On the ordinals, $j$ is increasing.  Thus, the second group of sentences tells us that $\kappa \leq a < j(\kappa)$, so $\crit j \leq \kappa$.  On the other hand, for each $\alpha < \kappa$, $V$ models the $\bL_{\kappa, \omega}$-sentence
$$\forall x \left(x E c_\alpha \to \bigvee_{\beta<\alpha} x = c_\beta\right)$$
We can use this to show that $c_\alpha = \alpha$ for every $\alpha < \kappa$ using induction.  Thus, $\crit j = \kappa$.

Now we claim that $b = j"\lambda$.  For each $\alpha \in \lambda$, we ensure that $j(\alpha) E b$, so $j"\lambda \subset b$.  Given $x \in b$, since $\cM$ omits $p$, we must have $x = j(\alpha)$ for some $\alpha < \lambda$.\hfill\dag\\

Note that we didn't directly use compactness for type omission in this proof.  However, in each case, the theory has a natural filtration by the appropriate partial order that is easily seen to be locally consistent while omitting the necessary type in $ZFC$.  For instance, in the case of $\kappa$ being $\lambda$-supercompact, for $s \in \cP_\kappa \lambda$, set $\alpha_s := otp(s)$.
\begin{eqnarray*}
T_s &:=& ED_{\bL_{\kappa, \kappa}}(V, \in, x)_{x \in (V_{\geq\kappa} \cup V_{\alpha_s})} \cup \{ c_i < c < c_\kappa \mid i < \alpha\} \cup \{c_i E d \wedge |d| < c_\kappa\mid i \in s\}\\
p_s(x) &:=& \{x E d \wedge x \neq c_i \mid i \in s\}
\end{eqnarray*}
Then $V$ is a model of this theory by interpreting every constant in the language by its index, $c$ as $\alpha_s$, and $d$ as $s$.  This gives a way to go directly between model-theoretic and elementary embedding characterizations.  It also shows that it is enough to omit a single\footnote{In the case of $n$-huge, recall that the omission of finitely many types can be coded by the omission of a single type.} type to obtain the $I$-$\kappa$-type omission for any number of types.

The ability to characterize cardinals at the level of huge and above shows that the addition of type omission to attempts to characterize large cardinals is a real necessity.  Measurable and strongly compact cardinals have known model-theoretic characaterizations without type omission, so one might wonder if type omission is necessary to characterize huge cardinals.  From the following theorem of Makowsky, we can deduce that it is necessary.

\begin{fact}[{\cite[Theorem 2]{m-vopcomp}}]\label{stavi-fact}
The following are equivalent:
\begin{enumerate}
	\item Every logic $\cL$ has a strong compactness cardinal; that is, for every logic $\cL$, there is a cardinal $\mu_\cL$ such that for any language $\tau$ and $\cL(\tau)$-theory $T$, if every $T_0 \in \cP_{\mu_\cL} T$ has a model then so does $T$.
	\item Vop\v{e}nka's Principle.
\end{enumerate}
\end{fact}

Thus, Vop\v{e}nka's Principle ``rallies at last to force a veritable G\"{o}tterdammerung" for compactness cardinals for logics\footnote{With apologies to Kanamori \cite[p. 324]{kanamori}.}.  Nonetheless, $\kappa$ being almost huge implies that $V_\kappa$ satisfies Vop\v{e}nka's Principle.  Thus, if $\kappa$ is the first huge cardinal, then $V_\kappa$ is a model of 
\begin{center}
Vop\v{e}nka's Principle + ``Every logic is compact, but there are no $\mu \leq \lambda$ such that $\bL_{\mu,\mu}$ is $[\lambda]^\mu$-$\mu$-compact for type omission"
\end{center}
Indeed, other approaches to model-theoretic characterizations of large cardinals focused solely on compactness or reflection principles have yet to characterize huge cardinals.

\section{Second-order logic and beyond}\label{second-order-sec}

We now turn to characterizations based on logics beyond (or orthogonal to) $\bL_{\infty, \infty}$.  In the spirit of Theorem \ref{framework-thm}, we can characterize compactness for omitting types in second-order logic with a similar theorem.

\begin{theorem} \label{second-thm}
Let $\kappa \leq \lambda$ and $I \subset \cP(\lambda)$ be $\kappa$-robust.  The following are equivalent:
\begin{enumerate}
	\item $\bL^2\cup\bL_{\kappa, \omega}$ is $I$-$\kappa$-compact for type omission.
	\item $\bL^2_{\kappa, \kappa}$ is $I$-$\kappa$-compact for type omission.
	\item \label{413} For every $\alpha > \lambda$, there is some $j:V_\alpha \to V_\beta$ such that $\crit j = \kappa$, $j(\kappa)>\lambda$, and $j"\lambda \in j(I)$.
	\item For every $\alpha > \lambda$, there is some $j:V \to \cM$ such that $\crit j = \kappa$, $j"\lambda \in j(I)$, and $V_{j(\alpha)} \subset \cM$.
\end{enumerate}
Moreover, the first $\mu$ such that $\bL^2$ is $I$-$\kappa$-compact for type omission is the first $\mu$ that satisfies (\ref{413}) except with $\crit j = \mu$.

\end{theorem}

{\bf Proof:} (4) implies (3) and (2) implies (1) are immediate.  We show that (1) implies (3) implies (4) implies (2).\\

For (1) implies (3), fix $\alpha > \lambda$ and consider the $\bL^2 \cup \bL_{\kappa, \omega}$-theory and type
\begin{eqnarray*}
T &=& ED_{\bL_{\kappa, \omega}}(V_\alpha, \in, x)_{x\in V_\alpha} \cup \{c_i < c < c_\kappa \mid i < \kappa\} \cup \{\Phi\}\\
p(x) &=&  \{ x E d \wedge x \neq c_i \mid i < \lambda\}
\end{eqnarray*}
where $d$ is a new constant symbol and $\Phi$ is Magidor's sentence from Fact \ref{magidor-phi-fact}.  Then, we can filtrate this theory as 
\begin{eqnarray*}
T_s &=& ED_{\bL_{\kappa, \omega}}(V_\alpha, \in, x)_{x\in V_{\sup s \cup [\kappa,\alpha)}} \cup \{c_i < c < c_\kappa \mid i < \sup s\} \cup \{\Phi\}\\
p_s(x) &=&  \{ x E d \wedge x \neq c_i \mid i \in s\}
\end{eqnarray*}
For each $s \in I$, we have that the natural expansion $(V_\alpha, \in, x, s)_{x \in V_{\sup s \cup [\kappa, \alpha)}}$ models $T_s$ and omits $p_s$.  Thus, our compactness principle tells us there is a model of $T$ omitting $p$, which, after taking the transitive collapse, gives the desired $j:V_\alpha \to V_\beta$.\\

For (3) implies (4), fix $\alpha \geq \lambda$ and let $\alpha'$ be the next strong limit cardinal above $\alpha$.  Then there is $j:V_{\alpha'}\to V_\beta$ with $\crit j =\kappa$ and $j"\lambda \in j(I)$.  Then derive the extender $E$ of length $\beth_{j(\alpha)}$ to capture this embedding.  Forming the extender power of $V$ and taking the transitive collapse, we get $j_E:V \to \cM_E$ with the desired properties.\\

For (4) implies (2), let $\bar{T} = \{T_s \mid s \in I\}$ be an increasing filtration of the $\bL^2_{\kappa, \kappa}$-theory $T$ that respects the index and $\{p^a(x) \mid a \in A\}$ be a collection of types indexed as $p^a(x) = \{\phi^a_i(x) \mid i < \lambda\}$ such that there are a club of $s \in I$ with a model $M_s$ that models $T_s$ and omits each $p^a_s$.  Fix strong limit $\alpha \geq \lambda$ to be greater than the rank of these models, their power sets, and the function $f$ that takes each of these $s$ to $M_s$; form $j:V \to \cM$ with $\crit j = \kappa$, $j"\lambda \in j(I) \cap \cM$, and $V_{j(\alpha)} \subset \cM$.  Since the domain of $f$ contains a club, it includes $j"\lambda$.  Set $M_*:= j(f)(j"\lambda)$.  By the elementarity of $j$, inside of $\cM$ we have that $M_* \vDash j(\bar{T})_{j"\lambda}$ and, for each $a \in A$, $M_*$ omits $j(p^a)_{j"\lambda} = \{j(\phi^a_i) \mid i < \lambda\} = j"p^a$.  Since $V_{j(\alpha)} \subset \cM$ and $\rank M_* < j(\alpha)$, $\cM$ is correct about this satisfaction.  Finally, $j"T \subset j(\bar{T})_{j"\lambda}$ because the filtration respects the index.  Thus, after a suitable renaming, we have found a model of $T$ omitting $\{p^a(x) \mid a \in A\}$.\hfill \dag\\

To aid in the discussion of the implication of this theorem, we introduce the following ad hoc naming convention for large cardinal properties.

\begin{defin}\label{ad-hoc-def}
Suppose a large cardinal property $P$ is characterized by being an $I$-$\kappa$-compactness cardinal for $\mathbb{L}_{\kappa,\kappa}$.  Given a logic $\cL$, we say that $\kappa$ is \emph{$P$-for-$\cL$} iff $\cL$ is $I$-$\kappa$-compact for type omission.
\end{defin}

For instance, Corollary \ref{main-cor}.(\ref{h}) characterizes huge as the existence of a $\lambda > \kappa$ such that $\bL_{\kappa, \kappa}$ is $[\lambda]^\kappa$-$\kappa$-compact for type omission, so saying that $\kappa$ is huge-for-$\bL^2_{\kappa, \kappa}$ means that there is a $\lambda > \kappa$ so $\bL^2_{\kappa, \kappa}$ is $[\lambda]^\kappa$-$\kappa$-compact for type omission.

Comparing Theorems \ref{framework-thm} and \ref{second-thm}, a large difference is that the first-order characterizations are witnessed by a single embedding, while the second-order characterizations require class many embeddings.  The reason for this is that a single model $\cM$ can be right about $\bL_{\kappa, \kappa}$ everywhere, but cannot be right about $\bL^2$ everywhere; otherwise, it would compute the power set of every set correctly and would be $V$.  Similarly, the type omitting compactness does \emph{not} hold for definable class theories for second-order as it does for first.  If it did, one could easily derive an nontrivial embedding $j:V \to V$.  

The first consequence of Theorem \ref{second-thm} regards the identity crisis.  In the language of Definition \ref{ad-hoc-def}, Magidor has shown that extendible cardinals are exactly those that are strong compact-for-$\bL^2_{\kappa,\kappa}$ \cite[Theorem 4]{m-roleof} and additionally shown that the first strongly compact cardinal could be the first measurable or the first supercompact\cite{m-ident-crisis}.  This second result means that various compactness notions for $\bL_{\kappa, \kappa}$ have an imprecise relation to one another: chain compactnes could coincide with compactness, or there could be many chain compact cardinals below the first compactness cardinal.  Surprisingly, when moving to $\bL^2$, these notions coincide and the identity crisis disappears!

\begin{theorem}\label{ext-char-thm}
The following are equivalent.
\begin{enumerate}
	\item $\kappa$ is measurable-for-$\bL^2 \cup \bL_{\kappa,\omega}$.
	\item $\kappa$ is strongly compact-for-$\bL^2_{\kappa,\kappa}$.
	\item $\kappa$ is supercompact-for-$\bL^2_{\kappa,\kappa}$.
\end{enumerate}
In particular, all three of these statements characterize extendible cardinals.
\end{theorem}

Here we take `$\kappa$ is measurable-for-$\cL$' as in Fact \ref{meas-fact}.  That is, we don't incorporate any type omission; however, the type omission characterization holds as a result of the above.

{\bf Proof:}  Clearly, (3) implies (2) implies (1) using (for the first implication) the trivially omitted type $\{x \neq x \mid i < \lambda\}$.  

The condition Theorem \ref{second-thm}.(3) is clearly stronger than extendibility, so any compactness for $\bL^2_{\kappa, \kappa}$ (including chain compactness) gives extendibility.  In particular, $j"\kappa = \kappa \in j(\cP_\kappa \kappa)$.  So measurable-for-$\bL^2 \cup \bL_{\kappa, \omega}$ implies extendible.

Similarly, the definition of extendibility includes that $j(\kappa) > \alpha$.  In this case, $j"\lambda$ has size $\lambda \leq \alpha$, so $j"\lambda \in j(\cP_\kappa \lambda)$.  Thus extendibility implies supercompact-for-$\bL^2_\kappa, \kappa$.\hfill \dag\\

The key to these equivalences is that the condition about $j"\lambda$ in Theorem \ref{framework-thm}.(\ref{emb}) often had more to do with the closure of the target model (i.e., is $j"\lambda$ in $\cM$?), rather than the nature of the relationship between $j(I)$ and $j"\lambda$.  When we have extendible-like embeddings, $j"\lambda$ is always in the target model, so many of the type omitting compactness principles (or even just compactness principles) become trivial.

A possible explanation for the collapse of the identity crisis is that type omission in $\bL^n_{\kappa, \kappa}$ is expressible\footnote{This is immediate for $\phi$-type omission for fixed $\phi$, and any type omission can be coded as $\phi$-type omission in an expansion.} in $\bL^{n+1}_{\kappa, \kappa}$, which is again codeable in $\bL^n_{\kappa, \kappa}$.  Thus, one might expect no difference between strong compact- and supercompact-for-$\bL^2_{\kappa, \kappa}$.  However, this does not explain why measurability coincides with these notions,  and the below proposition shows that some notions of type-omitting compactness for $\bL^2$ are strictly stronger than extendibility (in consistency strength).

\begin{prop}\label{huge-second-prop} \
\begin{enumerate}
	\item $\kappa$ is huge at $\lambda$-for-$\bL^2_{\kappa, \kappa}$ iff for every $\alpha \geq \lambda$, there is $j:V_\alpha \to V_\beta$ such that $\crit j = \kappa$ and $j(\kappa) = \lambda$.
	\item If $\kappa$ is almost 2-huge at $\lambda_1, \lambda_2$, then there is a $\kappa$-complete, normal ultrafilter on $\kappa$ containing
	$$\{\alpha < \kappa \mid V_{\lambda_2} \vDash ``\alpha \text{ is huge-for-}\bL^2_{\kappa, \kappa}"\}$$
	\item If $\kappa$ is huge-for-$\bL^2_{\kappa, \kappa}$, then there is a $\kappa$-complete, normal ultrafilter on $\kappa$ containing
	$$\{ \alpha < \kappa \mid \alpha \text{ is huge}\}$$
\end{enumerate}
\end{prop}

{\bf Proof:} The first item is just a restatement of Theorem \ref{second-thm} with $I = [\lambda]^\kappa$.  

Suppose $\kappa$ is 2-huge at $\lambda_1$, $\lambda_2$ and $j:V \to \cM$ witnesses this.  Fixing $\alpha \in [\lambda_1, \lambda_2)$, $j\rest V_\alpha$ is an embedding from $V_\alpha$ to $V_\beta$ with $j(\kappa) = \lambda_1$ that is in $\left(V_{j(\lambda_2)}\right)^{\cM}$.  So
$$\left(V_{j(\lambda_2)}\right)^{\cM} \vDash ``\exists j_0:V_\alpha \to V_\beta\text{ such that }\crit j_0 = \kappa\text{ and }j_0(\kappa)=\lambda_1"$$
Recall that\comment{ (by arguments of \cite{srk})}, $V_{\lambda_2} = \left(V_{\lambda_2}\right)^{\cM} \prec \left(V_{j(\lambda_2)}\right)^{\cM}$.  Thus,
$$V_{\lambda_2} \vDash ``\exists j_0:V_\alpha \to V_\beta\text{ such that }\crit j_0 = \kappa\text{ and }j_0(\kappa)=\lambda_1"$$
Since $\alpha$ was arbitrary,
\begin{eqnarray*}
V_{\lambda_2}&\vDash& ``\forall \alpha \geq \lambda_1, \exists j_0:V_\alpha \to V_\beta\text{ such that }\crit j_0 = \kappa\text{ and }j_0(\kappa)=\lambda_1"\\
V_{\lambda_2} &\vDash& ``\kappa \text{ is huge at }\lambda_1\text{-for-}\bL^2_{\kappa, \kappa}"
\end{eqnarray*}
Thus, $\{\alpha < \kappa \mid V_{\lambda_2} \vDash ``\alpha \text{ is huge-for-}\bL^2_{\kappa, \kappa}"\}$ is in the normal ultrafilter on $\kappa$ derived from $j$.

Suppose $\kappa$ is huge at $\lambda$-for-$\bL^2_{\kappa, \kappa}$.  Picking $\alpha$ large enough and getting the corresponding $j:V_\alpha \to V_\beta$ with $j(\kappa) = \lambda$, we can derive a normal, $\kappa$-complete, fine ultrafilter $U$ on $[\lambda]^\kappa$.  Then $U \in V_\beta$, so $V_\beta \vDash ``\kappa$ is huge."  Thus, $\{\alpha < \kappa \mid \alpha \text{ is huge}\}$ is in the normal, $\kappa$-complete ultrafilter on $U$ generated from $j$.\hfill\dag\\

Similar results show that $n$-huge-for-$\bL^2_{\kappa, \kappa}$ lies strictly between $n$-huge and almost $n+1$-huge.  The preceding argument is due to Gabriel Goldberg, who also reports that he can show that huge-for-$\bL^2_{\kappa, \kappa}$ can be characterized in terms of hyperhugeness. \footnotei{Make a fact, check Gabe's argument.  Can we adapt it to get ``huge at $\lambda$-for-$\bL^{s, \Sigma_n}_{\kappa, \kappa}$ implies $C^{(n)+}$-hyperhuge?"} 

Recall that $\kappa$ is $\lambda$-hyperhuge iff there is $j:V \to \cM$ with $\crit j = \kappa$ and ${}^{j(\lambda)} \cM \subset \cM$ and $\kappa$ is hyperhuge iff it is $\lambda$-hyperhuge for every $\lambda$.  Hyperhuge cardinals have recently been shown to imply the existence of a minimal inner model of $V$ that can reach $V$ by set-forcing extensions by Usuba \cite{u-ddg}.  Goldberg proves that $\kappa$ being hyperhuge is equivalent to the existence of a $\kappa_0<\kappa$ such that $\kappa_0$ is huge at $\kappa$-for-$\bL^2_{\kappa, \kappa}$.  Additionally, $\kappa$ being $\lambda$-hyperhuge is equivalent to the existence of $\mu > \lambda$ and a normal, fine, $\kappa$-complete ultrafilter on 
$$[\mu]^\lambda_{*\kappa}:=\{ s \subset \mu \mid |s| = \lambda, |s \cap \kappa| \in \kappa, \otp(s\cap\lambda) < \kappa\}$$
, which is equivalent to $\bL_{\kappa, \kappa}$ being $[\mu]^\lambda_{*\kappa}$-$\kappa$-compact for type omission by Theorem \ref{framework-thm}.

Examining the proof of Theorem \ref{second-thm}, we see that a level-by-level characterization of, e.g., $\alpha$-extendibility is harder due to the tricky nature of the L\"{o}wenheim-Skolem number for second-order logics.  In first-order, the L\"{o}wenheim-Skolem number of $\bL_{\lambda, \kappa}$ for theories of size $\mu$ is $\left((\lambda+\mu)^{<\kappa}\right)^+$, which is also its L\"{o}wenheim-Skolem-Tarski number.  For second-order logic, $LS(\bL^2)$ (for sentences) is the supremum of all $\Pi_2$-definable ordinals (V\"{a}\"{a}n\"{a}nen \cite[Corollary 4.7]{v-sort-intro}) and $LST(\bL^2)$ is the first supercompact, if one exists \cite[Theorem 2]{m-roleof} .  However, weak compactness restricts the size of the theory, so admits a more local characterization.  Denote the L\"{o}wenheim-Skolem number of sentences of $\bL^2_{\kappa,\kappa}$ by $\ell^2_\kappa$.

\begin{theorem}\label{wc-second-thm}
The following are equivalent for $\kappa$.
\begin{enumerate}
	\item $\kappa$ is weakly compact-for-$\bL^2 \cup \bL_{\kappa, \omega}$.
	\item $\kappa$ is weakly compact-for-$\bL^2_{\kappa, \kappa}$.
	\item Given any $\kappa+1 \subset \cM \subset V_{\ell^2_\kappa}$ of size $\kappa$, there is a partial elementary embedding $j:V_{\ell^2_\kappa} \to V_\beta$ for some $\beta$ with $\dom j = \cM$ and $\crit j = \kappa$.
\end{enumerate}
\end{theorem}

{\bf Proof:}  Clearly, (2) implies (1).  We show (1) implies (3) implies (2).\\

Suppose $\kappa$ is weakly compact-for-$\bL^2 \cup \bL_{\kappa, \omega}$ and let $\kappa+1 \subset \cM \subset V_{\ell^2_\kappa}$.  Let $T$ be the $\bL^2$ theory consisting of
\begin{enumerate}
	\item the $\bL_{\kappa, \omega}$-elementary diagram of $\cM$ in $V_{\ell^2_\kappa}$;
	\item $c_i < c < c_\kappa$ for $i<\kappa$; and 
	\item Magidor's $\Phi$ from Fact \ref{magidor-phi-fact}.
\end{enumerate}
Then every $<\kappa$-sized subset of $T$ is satisfiable as witnessed by an expansion of $V_{\ell^2_\kappa}$.  By weak compactness, we get a model of $T$, which must be some $V_\beta$.  This induces a partial function $j:V_{\ell^2_\kappa} \to V_\beta$ with $\dom j = \cM$.  Moreover, the elements of $T$ make this a partial elementary embedding with $\crit j = \kappa$.\\

Suppose that $\kappa$ satisfies the embedding property.  Let $T = \{\phi_i \mid i < \kappa\}$ be a $\bL^2_{\kappa, \kappa}(\tau)$-theory that is $<\kappa$-satisfiable with $|\tau|\leq \kappa$.  Then, there is a function $f$ with domain $\kappa$ such that $f(\alpha)\vDash \{\phi_i \mid i < \alpha\}$ for every $\alpha < \kappa$; moreover, by the definition of $\ell^2_{\kappa}$, we can assume that $f(\alpha) \in V_{\ell^2_\kappa}$.  Let $\cM \subset V_{\ell^2_\kappa}$ contain all of this information and be of size $\kappa$.  Then, there is partial elementary $j:V_{\ell^2_\kappa} \to V_\beta$ with $\dom j = \cM$ and $\crit j = \kappa$.  In particular, we have that 
\begin{enumerate}
	\item $j(\kappa) > \kappa$;
	\item $V_\beta \vDash ``j(f)(\kappa) \vDash j"T"$ and $V_\beta$ is correct about this; and
	\item $j"T$ and $T$ are just renamings of the same theory.
\end{enumerate}
Thus, the suitably renamed $j(f)(\kappa)$ witnesses that $T$ is satisfiable.\hfill\dag\\

A key piece in translating weak compactness for second-order into an embedding characterization is the ability to axiomatize well-foundedness.  If we look at a fragment of $\bL^2_{\kappa, \kappa}$ that includes an expression of well-foundedness, then weak compactness for this fragment is characterizable in a similar way, replacing $\ell^2_\kappa$ with the L\"{o}wenheim-Skolem number of that fragment.  However, if the fragment cannot express well-foundedness, then this characterization is harder.

Similar results can be proved by restricting the size of the theories under consideration.  In the general scheme, the theory $T$ is allowed to be as large as one wants, as are the pieces $T_s$ of the filtration.  If one restricts these pieces to be of size $\leq \mu$ and wants to characterize $\bL^2_{\kappa, \kappa}$ being $I$-$\kappa$-compact for type omission, then it suffices to look at an embedding as in Theorem \ref{second-thm}.(2) for $\alpha$ equal to the L\"{o}wenheim-Skolem number of $\bL^2_{\kappa, \kappa}$ for $\mu$-sized theories.

For the characterizations of strong and its variants, we need the concept of a Henkin second-order structure that is full up to some rank.  Recall the notion of a Henkin model described in Section \ref{prelim-sec}.

\begin{defin}\label{full-to-def}
Let $M_* = (M, P, E)$ be a Henkin structure and $A$ a transitive set.  
\begin{enumerate}
	\item $M_*$ is \emph{full to $A$} iff every $X \in \cP(M) \cap A$ is represented in $P$; this means that there is $c_X \in P$ such that, for all $y \in M$,
	$$y \in X \iff M_* \vDash y E c_X$$
	\item $M_*$ is \emph{full up to rank $\alpha$} iff it is full to $V_\alpha$.
\end{enumerate}
\end{defin}

While a Henkin structure has a nonstandard interpretation of second-order quantifiers, other additions to the logic must be interpreted standardly.  In particular, the next theorem discusses Henkin models of $\bL^2(Q^{WF})$-theories; while any second-order assertions of well-foundedness--i.e., $``\forall X \exists y \forall z (y \in X \wedge z \in X \to y R z)"$--can be satisfied non-standardly, any $Q^{WF}$ assertions of well-foundedness--i.e., $Q^{WF}xy(xRY)$--must be correctly interpreted (and so $R$ is well-founded).

\begin{theorem}\label{strong-char-thm}
The following are equivalent for $\kappa \leq \lambda$.
\begin{enumerate}
	\item $\kappa$ is $\lambda$-strong.
	\item \label{2strong} If $T \subset \bL^2_{\kappa, \omega}(Q^{WF})(\tau)$ is a theory that can be written as an increasing union $T = \cup_{\alpha < \kappa} T_\alpha$ such that every $T_\alpha$ has a (full) model, then $T$ has a Henkin model whose universe is an ordinal and is full up to rank $\lambda$.
	\item \label{3strong} Same as (\ref{2strong}), but there is also a type $p=\{\phi_i(x) \mid i <\kappa\}$ such that $T_\alpha$ has a (full) model omitting $p_{\alpha}$, and the resulting model omits $p$.
\end{enumerate}
\end{theorem}

Note that we add the condition on the universe of the model in (\ref{2strong}) to remove the possibility that the ``full up to rank $\lambda$" condition is vacuous; if the universe of $M$ just consists of elements of rank bigger than $\lambda$, then $M$ is trivially full up to rank $\lambda$.

{\bf Proof:} First, suppose that $\kappa$ is $\lambda$-strong and let $T$ be a theory and $p$ a type as in (\ref{3strong}).  We produce a model of $T$ in the standard way: let $f$ be a function with domain $\kappa$ such that $f(\alpha)$ is a model of $T_\alpha$.  WLOG, $f(s)$ is a full Henkin structure.  Then, in $\cM$, $j(f)(\kappa)$ is a model of (a theory containing) $j"T$ and $\cM \vDash ``j(f)(\kappa)$ is a full Henkin structure that omits $j(p)_\kappa = j"p"$.  $\cM$ is incorrect about second-order satisfaction above rank $\lambda$; however, since $V_\lambda \subset \cM$, it is correct about second-order satisfaction up to rank $\lambda$.\\

Second, suppose we have compactness.  Then we wish to build an embedding witnessing strength.  By the normal arguments, e.g. \cite[Section 26]{kanamori} or see Proposition \ref{strong-char-prop}, it is enough to derive a $(\kappa, \beth_\lambda)$-extender from an embedding $j: V_{\kappa+2} \to \cM$ with $\crit j = \kappa$, $V_\lambda \subset \cM$, and $\cM$ well-founded.  We can find such a model by considering the theory
$$ED_{\bL_{\kappa, \omega}}(V_{\kappa+2}, \in, x)_{x \in V_{\kappa+2}} \cup \{c_i < c < c_\kappa\mid i < \kappa\} \cup \{\Phi\} \cup\{Q^{WF}xy (x E y)\}$$
This can be written as an increasing $\kappa$-length union of satisfiable theories in the standard way and any model leads to, after taking transitive collapse, the necessary $j:V_{\kappa+2} \to \cM$.
\hfill \dag\\

We could ask for a variation of $(\ref{2strong})$ that allows for arbitrary $\kappa$-satisfiable theories or, equivalently, theories indexed by some $\cP_\kappa \mu$.  This would be equivalent to a jointly $\lambda$-strong and $\mu$-strongly compact cardinal: there is a $j:V \to \cM$ such that $\crit j = \kappa$, $V_\lambda \subset \cM$, and there is $Y \in j(\cP_\kappa \mu)$ such that $j"\mu \subset Y$.

If we drop the $Q^{WF}$, then we can characterize a weakening of $\lambda$-strong.  In the following theorem and proof, we break the convention that $\cM$ always denotes some transitive model of a fragment of $ZFC$.  In particular, we allow it to be ill-founded.  For such models, $wfp(\cM)$ denotes the well-founded part of $\cM$.

\begin{theorem} \label{ns-strong-thm}
The following are equivalent for $\kappa \leq \lambda$.
\begin{enumerate}
	\item $\kappa$ is non-standardly $\lambda$-strong: there is an elementary embedding $j: V\to \cM$ \emph{with $\cM$ not necessarily transitive} such that $\crit j = \kappa$ and $V_\lambda \subset wfp(\cM)$.
	\item If $\tau$ is a language and $T \subset \bL^2_{\kappa, \omega}$ is a theory that can be written as an increasing, continuous union $T = \cup_{\alpha < \kappa} T_\alpha$ such that every $T_\alpha$ has a (full) model, then $T$ has a Henkin model whose universe is an ordinal and is full up to rank $\lambda$.
\end{enumerate}
\end{theorem}

{\bf Proof:} The proof is the same as Theorem \ref{strong-char-thm}, with the changes exactly that we no longer insist on being correct regarding statements about well-foundedness.\hfill \dag\\

An argument of Goldberg shows that the level-by-level notions of non-standardly $\lambda$-strong and $\lambda$-strong are not equivalent, but full non-standard strong is equivalent to strong.

\subsection{$C^{(n)}$ and sort logic}

Moving to sort logic, we can prove a metatheorem along the lines of Theorems \ref{framework-thm} and \ref{second-thm} by introducing the notion of a $C^{(n)}$-cardinal.  The $C^{(n)}$ variants of large cardinals were introduced by Bagaria \cite{b-cn}.  Briefly, set $C^{(n)} = \{ \alpha \in \ON \mid V_\alpha \prec_{\Sigma_n} V\}$, where $\prec_{\Sigma_n}$ is elementarity for $\Sigma_n$ formulas in the L\'{e}vy hierarcy (in the language of set theory).  For a large cardinal notion $P$ witnessed by a certain type of elementary embedding, $\kappa$ is $C^{(n)}$-$P$ iff there is an elementary embedding $j$ witnessing that $\kappa$ is $P$ and so $j(\kappa) \in C^{(n)}$.  Recently, Tsaprounis \cite[Corollary 3.5]{t-cn-ext} and Gitman and Hamkins \cite[Theorem 15]{gh-genvop} have independently shown that $C^{(n)}$-extendibility is equivalent to the \emph{a priori} stronger notion of $C^{(n)+}$-extendibility: $\kappa$ is $C^{(n)+}$-extendible iff for all $\alpha > \kappa$ in $C^{(n)}$, there $j:V_\alpha \to V_\beta$ with $\crit j = \kappa$ and $j(\kappa), \beta \in C^{(n)}$.  It is the notion of $C^{(n)+}$-extendibility that we will use.

For some large cardinal notions, there is no increase of strength from moving to the $C^{(n)}$-versions (measurable, strong \cite[Propositions 1.1 and 1.2]{b-cn}, strongly compact \cite[Theorem 3.6]{t-exten-char}), but several other notions give an increasing hierarchy of strength.   Recall the notions of sort logic described in Section \ref{prelim-sec}

\begin{theorem}\label{sort-thm}
Let $\kappa \leq \lambda$, $n < \omega$, and $I \subset \cP(\lambda)$ be $\kappa$-robust.  The following are equivalent:
\begin{enumerate}
	\item $\bL^{s, \Sigma_n}\cup\bL_{\kappa, \omega}$ is $I$-$\kappa$-compact for type omission.
	\item $\bL^{s, \Sigma_n}_{\kappa, \kappa}$ is $I$-$\kappa$-compact for type omission.
	\item For every $\alpha \geq \lambda$ in $C^{(n)}$, there is some $j:V_\alpha \to V_\beta$ such that $\crit j = \kappa$, $j"\lambda \in j(I)$, and $\beta \in C^{(n)}$.
	\item For every $\alpha \geq \lambda$ in $C^{(n)}$, there is some $j:V \to \cM$ such that $\crit j = \kappa$, $j"\lambda \in j(I) \cap \cM$,$V_{j(\alpha)} \subset \cM$, and $j(\alpha) \in C^{(n)}$.
\end{enumerate}
\end{theorem}

The proof of Theorem \ref{sort-thm} follows the structure of Theorems \ref{framework-thm} and \ref{second-thm}.  To make the necessary changes, we introduce the following notion and lemma.  Given a $\Sigma_n$ formula $\phi(\bx)$ (in the L\'{e}vy hierarchy), let $\phi^\sim(\bx) \in \bL^{s, \Sigma_n}$ be the same formula where unbounded quantifiers are replaced with the corresponding sort quantifiers.  This allows us to characterize $C^{(n)}$ as follows.

\begin{lemma}
Let $\alpha$ be an ordinal.  $\alpha \in C^{(n)}$ iff $V_\alpha$ models
$$\{\forall \bx \left(\phi(\bx) \leftrightarrow \phi^\sim(\bx)\right) \mid \phi \text{ is }\Sigma_n\}$$
\end{lemma}

{\bf Proof:} For $\ba \in V_\alpha$, we always have $\phi(\ba)$ holds in $V$ iff $V_\alpha \vDash \phi^\sim(\ba)$.  The above theory makes this equivalent to $V_\alpha \vDash \phi(\ba)$. \hfill \dag\\

{\bf Proof of \ref{sort-thm}:}  We sketch the proof and highlight the changes from the proof of Theorem \ref{second-thm}.\\

Given the compactness, we prove (3) by considering the theory and type
\begin{eqnarray*}
T &=& ED_{\bL_{\kappa, \omega}}(V_\alpha, \in, x)_{x\in V_\alpha} \cup \{c_i < c < c_\kappa \mid i < \kappa\} \cup \{\Phi\}\\
& & \cup \{\forall \bx \left(\phi(\bx) \leftrightarrow \phi^\sim(\bx)\right) \mid \phi \text{ is }\Sigma_n\}\\
p(x) &=&  \{ x E d \wedge x \neq c_i \mid i < \lambda\}
\end{eqnarray*}
We filtrate this according to $I$ in the standard way and use expansion of $V_\alpha$ to provide witness models; here it is crucial that we started with $\alpha \in C^{(n)}$.  The model of $T$ omitting $p$ gives the desired $j$.  We can adjust this proof to get a proof of (4) by finding strong limit $\alpha' >\alpha$, also in $C^{(n)}$, and relativizing the appropriate parts of the theory to ensure that $j(\alpha) \in C^{(n)}$.  Then, derive the extender $E$ from this model, and $j_E: V \to \cM_E$ that retains the desired properties.\\

Given (3) or (4), we prove the compactness by starting with a filtration $\bar{T} = \{T_s \mid s \in I\}$ of an $\bL^{s, \Sigma_n}_{\kappa, \kappa}$-theory and types $\left\{p^a(x) = \{\phi^a_i(x) \mid i < \lambda\} \mid a \in A\right\}$, find strong limit $\alpha \in C^{(n)}$ above the rank of these objects and the function $f$ that takes $s$ to the model of $T_s$ omitting each $p^a_s$.  $V_\alpha$ reflects these properties since $\alpha \in C^{(n)}$, so by elementarity the target model thinks that $j(f)(j"\lambda)$ models $j"T$ and omits $\{j"p^a(x) \mid a \in A\}$.  Since $V_\beta$ or $V_{j(\alpha)}$ are $\Sigma_n$-elementary in $V$, the target model is correct.\hfill \dag\\

Similar to second-order logic, the identity crisis disappears in sort logic and $C^{(n)}$-extendible cardinals witness a wide range of type omitting compactness.
\comment{First, we show that $C^{(n)+}$-extendible cardinals are equivalent to an apparent weakening: the requirement that $j(\kappa) \in C^{(n)}$ is redundant.
\begin{prop}\label{cn-ext-prop}
The following are equivalent:
\begin{enumerate}
	\item $\kappa$ is $C^{(n)+}$-extendible.
	\item \label{4112} For every $\alpha > \kappa$ in $C^{(n)}$, there is $j:V_{\alpha} \to V_\beta$ such that $\crit j = \kappa$ and $\beta \in C^{(n)}$.
\end{enumerate}
\end{prop}}
We use the following lemma which will also be useful when examining L\"{o}wenheim-Skolem-Tarski numbers.  This is similar to Magidor's characterization of supercompacts.

\begin{lemma}\label{ext-lem}
Let $\kappa$ be $C^{(n)}$-extendible.  Then for all $\alpha >\kappa$ in $C^{(n)}$ and $R \subset V_\alpha$, there are cofinally many $\gamma < \kappa$ such that there are $\bar{\alpha} <\kappa$ in $C^{(n)}$ and $S \subset V_{\bar{\alpha}}$ with elementary $j:(V_{\bar{\alpha}}, \in, S) \to (V_{\alpha}, \in, R)$, $\crit j = \gamma$, and $j(\gamma) = \kappa$.
\end{lemma}

{\bf Proof:}  Fix $\alpha \in C^{(n)}$ above $\kappa$, $R \subset V_\alpha$, and $\beta < \kappa$.  Find $\alpha' >\alpha$ in $C^{(n)}$.  By assumption, there is $j:V_{\alpha'} \to V_{\beta'}$ with $\crit j = \kappa$, $j(\kappa)>\alpha'$, and $\beta' \in C^{(n)}$.  Given a transitive model $\cM$ of a fragment of $ZFC$, write $C^{(n), \cM}$ for $\cM$'s version of $C^{(n)}$.  Since $\alpha, \alpha' \in C^{(n)}$, $\alpha \in C^{(n), V_{\alpha'}}$.  By elementarity, $j(\alpha) \in C^{(n), V_{\beta'}}$.  Thus, 
\begin{eqnarray*}
V_{\beta'} \vDash ``\exists \bar{\alpha} < j(\kappa)\text{ and } S \subset V_{\bar{\alpha}}, j_0:\left(V_{\bar{\alpha}}, \in, S\right) \to \left(V_{j(\alpha)}, \in, j(R)\right)\text{ such that }\\j_0(\crit j_0) = j(\kappa), \crit j_0>j(\beta), \text{ and }\bar{\alpha} \in C^{(n)}"
\end{eqnarray*}
This is witnessed by $j \rest V_\alpha$.  By elementarity,
\begin{eqnarray*}
V_{\alpha'} \vDash ``\exists \bar{\alpha} < \kappa\text{ and } S \subset V_{\bar{\alpha}}, j_0:\left(V_{\bar{\alpha}}, \in, S\right) \to \left(V_{\alpha}, \in, R\right)\text{ such that }\\j_0(\crit j_0) = \kappa, \crit j_0>\beta, \text{ and }\bar{\alpha} \in C^{(n)}"
\end{eqnarray*}
This is the desired result; note that it implies $\bar{\alpha} \in C^{(n)}$ because $\alpha' \in C^{(n)}$.\hfill \dag\\

\comment{{\bf Proof of \ref{cn-ext-prop}:}  Suppose (2) holds.  By Lemma \ref{ext-lem}, the class of $C^{(n)}$-ordinals are cofinal in $\kappa$.  Since this class is club, $\kappa \in C^{(n)}$.  Fix $\alpha > \kappa$ in $C^{(n)}$ and find $j$ as in (2); it suffices to show $j(\kappa) \in C^{(n)}$.  Since $\kappa, \alpha \in C^{(n)}$, we have $\kappa \in C^{(n), V_\alpha}$.  By elementarity, $j(\kappa) \in C^{(n),V_\beta}$.  Since $\beta \in C^{(n)}$, this gives $j(\kappa) \in C^{(n)}$.\hfill\dag}

\begin{prop} \label{sort-comp-prop}
The following are equivalent for every $n<\omega$.
\begin{enumerate}
	\item $\kappa$ is $C^{(n)}$-extendible.
	\item $\kappa$ is measurable-for-$\bL^{s, \Sigma_n}\cup\bL_{\kappa, \omega}$.
	\item $\kappa$ is strong compact-for-$\bL^{s, \Sigma_n}_{\kappa, \kappa}$.
	\item $\kappa$ is supercompact-for-$\bL^{s, \Sigma_n}_{\kappa, \kappa}$.
\end{enumerate}
\end{prop}

{\bf Proof:}  This follows a similar argument as Theorem \ref{ext-char-thm}, just requiring that the $\alpha$'s be in $C^{(n)}$.\hfill \dag\\

However, the notion of a huge-for-$\bL^{s, \Sigma_n}_{\kappa, \kappa}$ cardinal would be similarly stronger in consistency strength than the notion in Proposition \ref{sort-comp-prop}.

While the L\"{o}wenheim-Skolem-Tarski number for second order was determined by Magidor in \cite{m-roleof} and Magidor and V\"{a}\"{a}n\"{a}nen have explored the L\"{o}wenheim-Skolem-Tarski numbers of various fragments of $\bL^2$ in \cite{mv-lst}, the L\"{o}wenheim-Skolem-Tarski number of sort logic seems unknown.  We give a characterization of these cardinals in terms of a $C^{(n)+}$-version of Magidor's characterization of supercompacts.  We work with L\"{o}wenheim-Skolem-Tarski numbers for strictly first-order languages to avoid the technicalities around trying to develop a notion of elementary substructure for sort logic.  See Section \ref{second-elem-subsec} for a definition of elementary substructure in second-order logic.

\begin{theorem}
The following are equivalent for $\kappa$:
\begin{enumerate}
	\item The conclusion of Lemma \ref{ext-lem}.
	\item For all $\alpha < \kappa$, if $N$ is a structure in a strictly first-order language $\tau$ of size $<\kappa$, then there is $M \prec_{\bL_{\alpha, \alpha}} N$ of size $<\kappa$ such that $M$ and $N$ have the same $\bL^{s, \Sigma_n}_{\alpha, \alpha}$-theory.
\end{enumerate}
\end{theorem}

{\bf Proof:} (1) implies (2): Let $N$ be a $\tau$-structure with $|\tau|<\kappa$.  Find $\gamma < \kappa$ that is above $\alpha$ and $|\tau|$ and find $\alpha'>\kappa$ such that $N \in V_{\alpha'}$.  Code the structure $N$ into a relation $R \subset V_{\alpha'}$  By assumption, there is $\bar{\alpha} \in (\gamma, \kappa)$ in $C^{(n)}$ and $j:(V_{\bar{\alpha}}, \in, S) \to (V_\alpha, \in, R)$ with $\crit j > \gamma$.  Then $S$ codes a structure $M$ in $V_{\bar{\alpha}}$ that, by elementarity, models $Th_{\bL^{s, \Sigma_n}_{\alpha, \alpha}}(N)$.  Moreover, $j \rest M$ is $\bL_{\alpha, \alpha}$-elementary, so the range of $j \rest M$ is the desired model.\\

(2) implies (1): Fix $\alpha > \kappa$ in $C^{(n)}$, $R\subset V_{\alpha}$, and $\beta < \kappa$.  Apply the assumption to the structure $(V_\alpha, \in, R, \delta)_{\delta < \beta}$ to get $M \prec_{\bL^{\beta, \beta}} N$ with the same $\bL^{s, \Sigma_n}$-theory.  This includes Magidor's $\Phi$, so after taking the transitive collapse, we get $\prec_{\bL_{\beta, \beta}}$-elementary $j:(V_{\bar{\alpha}}, \in, S, \delta)_{\delta<\beta} \to (V_\alpha, \in, R, \delta)_{\delta < \beta}$.   The constants for the elements of $\beta$ force the critical point of $j$ above $\beta$.\hfill \dag\\

\subsection{Rank-into-rank} \label{rir-sec}

We now turn to the strongest large cardinal principles, the rank-into-rank embeddings.  For an excellent survey of these, see Dimonte \cite{d-rank-survey-pub} or the expanded Dimonte \cite{d-rank-survey}.  Following \cite{d-rank-survey-pub}, this section uses $\Sigma^1_{n}(\bL_{\kappa,\omega})$ to denote the fragment of infinitary second-order logic $\bL^2_{\kappa, \omega}$ consisting of the formulas that are $\Sigma_n$ in their second-order quantifiers.

We can cast $I1$, $I2$, and $I3$ in a uniform way by saying, for $n<\omega$, $I2_n(\kappa, \delta)$ is the assertion
\begin{center}
There is $j:V_\delta \to V_\delta$, $\kappa = \crit j$, $\delta$ is the supremum of the critical sequence, and $j$ is $\Sigma^1_{2n}(\bL_{\kappa, \omega})$-elementary.
\end{center}
Laver \cite[Theorem 2.3]{l-imp} proved\footnote{Laver's paper attributes this result to Martin without citation, but other sources attribute it to Laver.} that $\Sigma^1_{2n+1}$-elementarity of such a $j$ implies its $\Sigma^1_{2n+2}$-elementarity, so it suffices to consider the even levels.  Then $I3$ is $I2_0$, $I2$ is $I2_1$\footnote{Which is in turn equivalent to being elementary about statements of well-foundedness \cite[Lemma 6.13]{d-rank-survey}.}, and $I1$ is $I2_{<\omega}$.

Note that second-order elementary embeddings should be understood in the context of Section \ref{second-elem-subsec}.  Given $j:V_\delta \to V_\delta$, we can naturally extend this to $j^+:\cP(V_\delta) \to \cP(V_\delta)$ by $j^+(R) = \bigcup_{\alpha <\delta}j(R \cap V_\alpha)$. Then set $A^{V_\delta} = j^+(A)$.  Note that $j^+$ is the only possible extension of $j$ to $V_{\delta+1}$ that could be elementary, and its $\Sigma^1_0$-elementarity follows from its first-order elementarity.

These principles have natural characterizations in terms of extendibility criteria.  Recall that $\kappa$ is weakly compact iff every $\kappa$-sized structure has a proper $\bL_{\kappa, \kappa}$-elementary extension \cite[Theorem 4.5]{kanamori} and $\kappa$ is measurable iff every $\geq \kappa$-sized structure has a proper $\bL_{\kappa, \kappa}$-elementary extension.

\begin{prop}
For each $n\leq\omega$ and $\delta = \beth_\delta$, the following are equivalent.
\begin{enumerate}
	\item There is $\kappa < \delta$ such that $I2_n(\kappa, \delta)$.
	\item Every $\delta$-sized structure has a non-identity $\Sigma^1_{2n}(\bL_{\omega, \omega})$-elementary embedding into itself.
\end{enumerate}
\end{prop}

{\bf Proof:} Suppose (1) holds, and let $j:V_\delta\to \delta$ with $\crit j =\kappa$ witness.  If $M$ is of size $\delta$, then we can code it as a structure with universe $\delta$.  Then $j\rest \delta$ is the desired embedding and it inherits the $\Sigma^1_{2n}(\bL_{\omega, \omega})$-elementarity of $j$.

Suppose (2) holds.  Then apply it to the structure $(V_\delta, \in)$, which has size $\beth_\delta=\delta$; this gives $\Sigma^1_{2n}(\bL_{\omega,\omega})$-elementary $j:V_\delta\to V_\delta$.  $j$ must have a critical point below $\delta$, call it $\kappa$.  As in Section \ref{prelim-sec}, this strengthens the elementarity of $j$ to $\Sigma_{2n}^1(\bL_{\kappa, \omega})$.  Since $j:V_\delta \to V_\delta$, we must have $\delta \geq \sup_n j^n(\kappa)$.  By Kunen's inconsistency, we must have $\delta \leq \sup_n j^n(\kappa)$.  Since $\delta$ is limit, this tells us that $\delta = \sup_n j^n(\kappa)$, and $I2_n(\kappa, \delta)$ holds.\hfill \dag\\

$I3$ and $I2$ also have standard characterizations in terms of coherent $\omega$-sequences of normal ultrafilters.  This allows us to prove the following type omitting compactness from them.  Unfortunately, this does not give an equivalence, but does allow us to sandwich these properties between two compactness for omitting types statements.  Recall that any type can be trivially extended to an equivalent, larger one by adding instances of ``$x \neq x$."

\begin{theorem}\
\begin{enumerate}
	\item If $I2_0(\kappa, \delta)$, then for any theory $T= \cup_{s \in \cP_\kappa \kappa}T_s \subset \bL_{\kappa,\omega}(\tau)$ and set of types $\{p^\beta = \{\phi^\beta_i(x) \mid i < \kappa_{n_\beta+1}\} \mid \beta < \mu\}$, if there are club many $s \in \cP_\kappa \kappa$ such that $T_{s}$ has a model $M_s$ with the property
\begin{center}
for $\beta < \mu$, $\left\{t \in [\kappa_{n_\beta+1}]^{\kappa_{n_\beta}} \mid M_{t \cap \kappa} \text{ omits } p^\beta_t = \{\phi^\beta_i(x) \mid i \in  t\}\right\}$ contains a club,
\end{center}
then $T$ has a model omitting $\{p^\beta \mid \beta < \mu\}$.
	\item $I2_1(\kappa, \delta)$ implies the same for the logic $\bL_{\kappa, \omega}(Q^{WF})$.
\end{enumerate}
\end{theorem}

We will use the ultrafilter characterizations of these cardinals, in part because $I3$ doesn't have a characterization in terms of $j:V \to \cM$ and we don't want to restrict to models of size $\leq\delta$.

{\bf Proof:} We prove the second item.  The first follows by the same argument, just removing the mentions of well-foundedness.

$I2_1(\kappa,\delta)$ is equivalent to the existence of $\kappa$-complete, normal ultrafilters $U_n$ on $[\kappa_{n+1}]^{\kappa_n}$ (where $\{\kappa_n\mid n < \omega\}$ is the critical sequence of the witnessing embedding) such that
\begin{enumerate}
	\item \emph{coherence:} For any $X \subset [\kappa_{n+1}]^{\kappa_n}$ and $m > n$,
	$$X \in U_n \iff \{s \in [\kappa_{m+1}]^{\kappa_m} \mid s \cap \kappa_{n+1} \in X\} \in U_m$$
	\item \emph{well-founded:} For any $\{n_i < \omega \mid i <\omega\}$ and $\{X_i \in U_{n_i}\mid i < \omega\}$, there is $s \subset \delta$ such that, for all $i < \omega$, $s \cap \kappa_{n_i+1} \in X_i$.
\end{enumerate}
Let $M_s$ be the desired models.  Since they exists for club many $s$ and $U_0$ contains this club, we can form the direct limit of the ultrapowers as standard: $M^*_n = \prod_{s \in [\kappa_{n+1}]^{\kappa_n}} M_{s \cap \kappa}/U$ and there is a coherent system of $\bL_{\kappa, \kappa}$-embeddings $f_{n, m}: M^*_n \to M^*_m$ that takes $[f]_{U_n}$ to $[s \mapsto f(s \cap \kappa_{n+1})]_{U_{n+1}}$.  Then $M^*$ is the direct limit of these models.  Standard arguments (see Proposition \ref{wf-los-prop} for a more general case) show that \L o\'{s}' Theorem holds for formulas of $\bL_{\kappa, \omega}(Q^{WF})$.  This guarantees that $M^* \vDash T$.  To show the type omission, let $\beta < \mu$ and $[n, f]_{\bar{U}} \in M^*$.  Setting $n' = \max{n, n_\beta}$, this implies
$$X:=\left\{t \in [\kappa_{n'+1}]^{\kappa_{n'}} \mid M_{t \cap \kappa} \text{ omits } p^\beta_{t\cap \kappa_{n_\beta+1}}\right\} \in U_{n'}$$
  Define a function $h$ on $X$ by setting, for $t \in X$, $h(t) \in t$ such that $M_{t \cap \kappa} \vDash \neg \phi^\beta_{h(t)}\left(f(t\cap \kappa_{n+1})\right)$; this is possible exactly because of the type omission.  Then, $h$ is regressive on a $U_{n'}$-large set, so there is $\alpha_0 < \kappa_{n'+1}$ such that
$$Y:=\left\{t \in [\kappa_{n'+1}]^{\kappa_{n'}} \mid M_{t \cap \kappa} \vDash \neg \phi^\beta_{\alpha_0}\left(f(t\cap\kappa_{n+1})\right)\right\} \in U_{n'}$$
Then, by \L o\'{s}' Theorem, $[n,f]_{\bar{U}} = [n', f(s \cap \kappa_{n+1})]_{\bar{U}}$ omits $p^\beta$ in $M^*$ as 
$$M^* \vDash  \neg \phi^\beta_{\alpha_0}\left([n, f]_{\bar{U}}\right)$$
Since $[n, f]_{\bar{U}}$ and $\beta$ we arbitrary, were are done.\hfill\dag\\

We can also isolate a type omitting compactness stronger than these rank-into-rank axioms.  Note that, unlike previous theorems, the types in the following don't shrink in the hypothesis.

\begin{theorem} Fix a cardinal $\delta=\beth_\delta$.
\begin{enumerate}
	\item Suppose we have the following for some $\kappa$:
	\begin{center}
	For any $\bL_{\omega, \omega}(\tau)$-theory $T$ with a filtration $\{T_\alpha \mid \alpha < \kappa\}$ and $|\tau| = \delta$ and any $\bL_{\omega, \omega}(\tau)$-type $p(x)$, if every $T_\alpha$ has a model omitting $p$, then $T$ has a model omitting $p$.
	\end{center}
	Then there is $\kappa_0 \leq \kappa$ such that $I2_0(\kappa_0, \delta)$.
	\item The above implies $I2_n(\kappa, \delta)$ after replacing the logic with $\Sigma^1_{2n}(\bL_{\omega, \omega})$.
\end{enumerate}
\end{theorem}
{\bf Proof:} Fix a bijection $f:\delta \to V_\delta$. Consider the theory and types
\begin{eqnarray*}
T&=& ED_{\bL_{\omega, \omega}}(V_\delta, \in, x)_{x \in V_\delta} \cup \{c_i < c < c_\kappa\mid  i < \kappa\}\\
& & \cup \{d_i E d_j \mid i, j < \delta, f(i) \in f(j)\}\cup \{\neg(d_i E d_j) \mid i, j < \delta, f(i) \not\in f(j)\}\\
p(x) &=& \{x \neq d_i \mid i < \delta\}
\end{eqnarray*}
We claim that a model of $T$ omitting $p$ will give a bijection as required for $I2_0(\kappa_0,\delta)$: the first two parts of $T$ will give an elementary $j:V_\delta \to \cM$, where we don't yet know that $\cM$ is transitive.  The third and fourth parts of $T$ ensures that $\cM$ has a (different) copy of $V_\delta$ in it given by $x \in V_\delta \mapsto d^\cM_{f^{-1}(x)} \in \cM$.  Omitting $p$ means everything in $\cM$ is in this second copy of $V_\delta$, so $\cM \cong V_\delta$.  Thus, we have elementary $j:V_\delta \to V_\delta$ that is not the surjective because its range does not contain $c^\cM$.  If $\kappa_0 := \crit j$, then this witnesses that $I2_0(\kappa_0, \delta)$ holds; note the second part of $T$ ensures that $\crit j \leq \kappa$.\\

To find such a model, we apply our compactness principle.  If $\alpha < \kappa$, let $T_\alpha$ be $T$ cutting out the constants $c_x$ for elements with rank in the interval $[\alpha, \kappa)$ as in the discussion following Theorem \ref{elem-comp-thm}, that is,
\begin{eqnarray*}
T_\alpha&=& ED_{\bL_{\omega, \omega}}(V_\delta, \in, x)_{x \in V_\alpha \cup V_{\geq \kappa}} \cup \{c_i < c < c_\kappa\mid  i < \alpha\}\\
& & \cup \{d_i E d_j \mid i, j < \delta, f(i) \in f(j)\}\cup \{\neg(d_i E d_j) \mid i, j < \delta, f(i) \not\in f(j)\}
\end{eqnarray*}
We can find a model of $T_\alpha$ omitting $p$ by taking the identity from $V_\delta$ to itself and interpreting $c$ as $\alpha$.\hfill\dag\\

\subsection{Elementary substructure for $\bL^2$} \label{second-elem-subsec}

In contrast with first-order logics\footnote{Here meant as `sublogics of $\bL_{\infty,\infty}$.'}, the notion of elementary substructure in \emph{second-order} logic has not been well-studied.  As evidence to this, in V\"{a}\"{a}n\"{a}nen's paper on second-order logic and set theory \cite{v-2nd-set}, the word `structure' appears 210 times, but `substructure,' `sub-structure,' or even `sub structure' never appear.  Similarly, Shapiro's book \cite{s-2nd-book} never discusses the matter\footnote{The notion of ``elementary substructure" does appear here, but always in reference to its first-order version.}.  A guess at the cause for this is that second-order logic is often employed to find categorical theories, whereas first-order logic attempts only to axiomatize classes of structures.  In such contexts, there is no reason to talk about one model's relation with others.

A first-attempt at second-order elementary substructure would be to work in analogy with first-order and say $M$ is a $\bL_2$-elementary substructure of $N$ iff every formula holds of parameters from $M$ in $M$ iff it holds in $N$.  However, this notion doesn't allow for any proper elementary extensions as there are definable sets that must grow.  Concretely, the formula $\phi(X) := ``\forall x \left(x \in X\right)"$ must be satisfied by the entire universe, so no extension of $M$ can think $\phi(M)$ holds.  A more modest generalization is used by Magidor and V\"{a}\"{a}n\"{a}nen in \cite[Between Definitions 3 and 4]{mv-lst}.  They say that $M$ is an $\bL^2$-elementary substructure of $N$ iff the above property holds \emph{restricted} to formulas with only first-order free variables.  This works (in the sense that proper extensions can exist), and they observe that Magidor's theorem on the L\"{o}wenheim-Skolem-Tarski number of second-order extends to this notion of elementarity.  However, it seems lacking as there's no discussion of free variables in second order.

We give a definition of elementary substructure in second-order logic in Definition \ref{2nd-elem-def} below.  An important point is that this definition includes a stronger notion of substructure in the second-order context (see Definition \ref{2nd-sub-def}.(3)).  Before giving the formal definitions, we give a motivation.  Any second-order structure $M$ can be viewed as a Henkin structure $M^+ = (M, P[M], E)$.  Then every second-order statement about $M$ can be transferred to a first-order statement about $M^+$, and the fact that $M^+$ is a full structure (isomorphic to $(M, \cP(M), \in )$) can be captured by a single second statement $\Psi$ asserting every subset of $M$ is represented by a member of $E$.

Now that we have moved to a more familiar first-order setting, we can ask what relation between the original structures $M$ and $N$ characterizes when $M^+$ is an elementary substructure of $N^+$.  The key is that, given $s \subset M$ and $m_s \in P[M]$ representing it, $N^+$ thinks the same facts about $m_s$ as $M$ does when the parameters come from $M$, but $N$ also thinks new things about $m_s$.  Crucially, there might be elements of $N-M$ that $N$ thinks are in $m_s$.  Then, setting $s^N := \{n \in N \mid N \vDash n E m_s\}$, $N$ thinks all the same facts about $s^N$ that $M$ does about $s$.  This notion of extending subsets is key to defining second-order elementary substructure.  This leads to the following definitions\footnote{We focus on $\prec_{\bL^2}$ for ease, but these definitions could easily be changed to accommodate $\prec_{\bL^2_{\alpha, \alpha}}$ for $\alpha < \kappa$.}:

\begin{defin}\label{2nd-sub-def}\
\begin{enumerate}
	\item Set $\cP^\omega(X) = \cup_{n < \omega} \cP(X^n)$.
	\item For $X\subset Y$, an \emph{extension function} is a function $f:\cP^\omega(X) \to \cP^\omega(Y)$ such that, for all $s \subset X^n$, $f(s) \cap X^n = s$ and $f(s) \subset Y^n$ and $f$ fixes every finite set.
	\item Given $\tau$-structures $M$ and $N$, we say $M$ is a \emph{second-order substructure} of $N$, written $M \subset_2 N$ iff $|M| \subset |N|$ and there is an extension function $s \mapsto s^N$ such that for every atomic $\phi(x_1, \dots, x_n; X_1, \dots, X_{n'})$, $m_i \in M$, and $s_i \subset M^{n_i}$, we have
$$M \vDash \phi(m_1, \dots, m_n; s_1, \dots, s_{n'}) \iff N \vDash \phi(m_1, \dots, m_n; s_1^N, \dots, s_{n'}^N)$$
\end{enumerate}
\end{defin}

We turn substructure into elementary substructure by letting $\phi$ range over all formulas.

\begin{defin}\label{2nd-elem-def}
Given $\tau$-structures $M$ and $N$, we say $M$ is an $\mathbb{L}^2$-elementary substructure of $N$, written $M \prec_{\mathbb{L}^2} N$ iff $M \subset_2 N$ and there is an extension function $s \mapsto s^N$ witnessing this such that for every $\phi(x_1, \dots, x_n; X_1, \dots, X_m) \in \mathbb{L}^2$, $m_i \in M$, and $s_i \subset M^{n_i}$, we have
$$M \vDash \phi(m_1, \dots, m_n; s_1, \dots, s_{n'}) \iff N \vDash \phi(m_1, \dots, m_n; s_1^N, \dots, s_{n'}^N)$$\end{defin}

Thus, $M \prec_{\bL^2} N$ comes with a choice of extensions for each $s \subset M$.  This avoids the issue with the ``first-attempt" notion of elementary substructure above.  Indeed, for any definable subset $A$ of $M$, $A^N$ must be the collection of elements satisfying that definition in $N$ if $M \prec_{\bL^2} N$.

Note that if $\tau$ is a strictly first-order language (as is often the case), then $M \subset N$ is equivalent to $M\subset_2 N$ for any collection of extensions.  This means that comments about substructures in first-order languages in the context of second-order can be given the normal interpretation.  Also, for such a language, $\prec_{\bL^2}$ will imply elementarity in the sense of Magidor and V\"{a}\"{a}n\"{a}nen.

For a general notion of $\leq$, a $\leq$-embedding is $f:M \to N$ such that $f$ is a $\tau$-isomorphism onto its range and $f(M) \leq N$ (a set theorist might prefer to call this model $f"M$).  Specializing to $\prec_{\bL^2}$, $f$ is a map on elements of $M$ and, given $s \subset M$, $f"s \subset f(M)$.  Then $f(M) \prec_{\bL^2} N$ implies there is an extension $(f"s)^N \subset N$ that satisfies the definition of $\prec_{\bL^2}$.  To avoid this unfortunate notation, we say that $f:M \to N$ is ${\bL^2}$-elementary means that $f$ is a map from $M \cup \cP^\omega(M)$ to $N \cup \cP^\omega(N)$ such that, for all $\phi(\bx, \bX)$, $\ba \in M$, and $\bs \in \cP^\omega (M)$,
$$M\vDash \phi(\ba, \bs) \iff N \vDash \phi(f(\ba), f(\bs))$$

One should be skeptical that this is the ``right" notion of $\prec_{\bL^2}$ as it adds a strange new condition about global choices of extensions of subsets.  In addition to the heuristic with Henkin models above, we offer two ``sanity" checks that this is the right notion.  

The first check is that Magidor's Theorem on the LST number of second order logic holds with this notion.  

\begin{cor}\label{2nd-lst-cor}
Let $\kappa$ be supercompact.  For any $\tau$ of size $<\kappa$ and $\alpha < \kappa$ and $\tau$-structure $N$, there is $M \prec_{\bL^2_{\alpha, \alpha}} N$ of size $<\kappa$.
\end{cor}

This follows from the argument given in our heuristic: take a model $N$, turn it into a full Henkin structure $N^+$ with Skolem functions, and apply Magidor's version along with the sentence $\forall X \exists x \left(X \text{`}=\text{'} x\right)$.  This quotational equality is an abbreviation of the statement that $X$ and $x$ have the same elements.

The second check is that the natural notion of an elementary diagram characterizes ${\bL^2}$-elementary embedability.  Given a $\tau$-structure $M$, we define it's $\bL^2$ elementary diagram by first adding a first-order constant $c_a$ for each $a \in M$ and a second-order constant $d_s$ for each $s \in \cP^\omega(M)$ and setting
$$ED_{\bL^2}(M) = \left\{ \phi(c_{a_1}, \dots, c_{a_n}, d_{s_1}, \dots, d_{s_k}) \mid M \vDash \phi(a_1, \dots, a_n, s_1, \dots, s_k) \right\}$$

\begin{prop}
Let $M, N$ be $\tau$-structures.  The following are equivalent.
\begin{enumerate}
	\item $N$ has an expansion that models $ED_{\bL^2}(M)$.
	\item There is a ${\bL^2}$-elementary $f:M \to N$.
\end{enumerate}
\end{prop}

{\bf Proof:} The proof follows as the first-order one.  Given $\bL^2$-elementary $f:M\to N$, expand $N$ by $c^N_a = f(a)$ and $d^N_s = f(s)$.  Similarly, if $N^* \vDash ED_{\bL^2}(M)$, define $f:M\to N$ by the same formula.\hfill\dag\\

Some of the basic results of first-order elementary substructure transfer, while others do not.  For instance, the Tarski-Vaught test goes through with the same proof, although a slightly modified statement.

\begin{prop}[Tarski-Vaught Test for $\prec_{\bL^2}$]
Given $M \subset_2 N$ (with a specified extension map $s \mapsto s^N$), we have that $M \prec_{\bL^2} N$ iff both of the following hold:
\begin{enumerate}
	\item $\forall m_i \in M, s_j \subset |M|^{n_j}$ and $\phi(x, \bx, \bX)$, if $N \vDash \exists x \phi(x, \bm, \bs^N)$, then there is $m \in M$ such that $N \vDash \phi(m, \bm, \bs^N)$; and
	\item $\forall m_i \in M, s_j \subset |M|^{n_j}$ and $\psi(X, \bx, \bX)$, if $N \vDash \exists X \phi(X, \bm, \bs^N)$, then there is $s \subset M^n$ such that $N \vDash \phi(s^N, \bm, \bs^N)$.
\end{enumerate}
\end{prop}

Unfortunately, $\prec_{\bL^2}$ does not have properties such as coherence and smoothness under unions of chains fail.  These are properties coming from the study of Abstract Elementary Classes, a general framework for nonelementary model theory (see Baldwin \cite{baldwinbook} for a survey):
\begin{itemize}
	\item {\bf Coherence:} A binary relation $\prec$ on $\tau$-structures is \emph{coherent} iff whenever $M_0 \prec M_2$ and $M_1\prec M_2$ and $M_0 \subset M_1$, then $M_0 \prec M_1$.
	\item {\bf Smoothness:} A  binary relation $\prec$ on $\tau$-structures is \emph{smoot under unions of chains} iff whenever $\{M_i \mid i < \alpha\}$ is a $\prec$-increasing chain with $\alpha$ limit, then $\bigcup_{i<\alpha} M_i$ is the $\prec$-least upper bound of the chain.
\end{itemize}

We can show $\prec_{\bL^2}$ fails these properties by using Silver's example, the bane of many a nonelementary model theorist's hope.  Recall Silver's example\footnote{Although typically given as a $PC$-class (see \cite[Example on p. 92]{k-lw1w}), we give it here as the class of models of a second-order sentence.}: the language consists of two sorts $P$ and $Q$ and the sentence $\phi_{Silver}$ is
$$\exists E \subset P \times Q \left( \forall x,y \in Q \left(x = y \iff \forall z \in P (z E x \iff z Ey )\right)\right)$$
In other words, a model $M$ of $\phi_{Silver}$ consists of two sets $P^M$ and $Q^M$ such that there is an extensional relation on $P^M$ and $Q^M$.  The existence of such a relation is equivalent to there being an injection from $P^M$ to $\cP(Q^M)$, which is equivalent to $|Q^M| \leq 2^{|P^M|}$.

If $\cF$ is any fragment of $\bL^2$ containing $\phi_{Silver}$, then $M \prec_{\cF} N$ requires that any extensional relation on $P^M \times Q^M$ can be extended to an extensional relation on $P^N \times Q^N$.  This is easy unless $|Q^M| = 2^{|P^M|}$ and $|Q^N| > 2^{|P^N-P^M|}$, but it fails in this case.  Thus, we can find counter-examples to coherence and smoothness under unions of chains.  Note that, since the language is first-order, $M \subset_2 N$ is equivalent to normal substructure.
\begin{itemize}
	\item {\bf Coherence:}  Set $M_0, M_1, M_2$ by $P^{M_0} = P^{M_1} = \omega$; $P^{M_2} =\omega_1$; $Q^{M_0} = \cP(\omega)$; and $Q^{M_1} = Q^{M_2} = \cP(\omega+1)$.  Then $M_0 \subset M_1$ and $M_0, M_1 \prec_{\bL^2} M_2$, but $M_0 \not\prec_{\bL^2} M_1$.
	\item {\bf Smoothness:} For $\alpha < 2^\omega$, set $M_\alpha$ by $P^{M_\alpha} = \omega$ and $Q^{M_\alpha} = \omega+\alpha$.  Then this is a $\prec_{\bL^2}$-increasing, continuous chain with union $M = \left( \omega, 2^\omega\right)$.  Set $N$ by $P^N = \omega$ and $Q^N = 2^\omega+\omega$.  Then $M_\alpha \prec_{\bL^2} N$ for each $\alpha < 2^\omega$, but $M \not\prec_{\bL^2} N$.
\end{itemize}

We can also define Skolem function for second-order in the same way as first-order.  To do so, we must add strictly second-order functions to the language, so that even if $\tau$ started out strictly first-order, its Skolemization $\tau_{sk}$ won't be.

In the context of first-order, the notions of elementary substructure and that of club sets are very closely intertwined.  For instance, given a structure $M$ and cardinal $\kappa < \|M\|$, 
$$\left\{s \in \cP_\kappa \lambda : s \text{ is the universe of an elementary substructure of }M\right\}$$
 is club, and, conversely, given a club $\cC \subset \cP_\kappa \lambda$, we can find a structure $M$ with universe $\lambda$ such that the above set is contained in $\cC$.  This connection is mediated by the fact that both notions can be characterized by being the closure sets of certain functions (Skolem functions, in the case of elementary substructure).
 
 With a definition for second-order elementary substructure in hand, we can define a notion of club, which we call superclubs.  Recall our focus on $\bL^2$.

\begin{defin}\
\begin{enumerate}
\item Fix $\mu < \kappa$ and let $F_i:[\lambda]^{<\omega} \times [\cP^\omega(\lambda)]^{<\omega} \to \cP^\omega(\lambda)$ for $i < \mu$.  Then
\begin{eqnarray*}
C(\bar{F}):=\{ s \in \cP_\kappa \lambda \mid \exists \text{ extension }f:\cP^\omega(s) \to \cP^\omega(\lambda) \text{ such that }\forall \ba \in s, \bx \subset s^{\bn}, i < \mu,\\
 \text{ we have } F_i\left(\ba, f(\bx)\right) = f(y) \text{ for some }y \in \cP^\omega(s)\}
\end{eqnarray*}
\item We call the collection $\cF$ of all sets containing some $C(\bar{F})$ the \emph{superclub filter on $\lambda$}.
\end{enumerate}
\end{defin}

We have called $\cF$ a filter without proving it is one.  This name is justified in Proposition \ref{sc-filter-prop}, although it might be non-proper.

Although we defined it with a combinatorial characterization in the spirit of Definition \ref{basic-def}, our interest in superclubs comes from the following model theoretic characterization.  Given a structure $M$ and $s \subset |M|$, set $M\rest s$ to be the $\tau$-substructure of $M$ with universe generated by $s$ (so is the closure of $s$ under the functions of $M$).

\begin{lemma}\label{sc-gen-lem}
The superclub filter is generated by sets of the form
$$D(M):=\{ s \in \cP_\kappa \lambda \mid M\rest s \prec_{\bL^2} M\}$$
for $M$ a $\tau$-structure with universe $\lambda$ and $|\tau|<\kappa$.
\end{lemma}

{\bf Proof:}  First, suppose we are given $\{F_\alpha\mid \alpha < \mu\}$.  Set $\tau =\{F_\alpha^{n; m_1, \dots, m_k} \mid \alpha <\mu; k, n, m_i < \omega\}$ to be a functional language so the domain of $F_\alpha^{n; m_1, \dots, m_k}$ is $M^n \times \cP(M^{m_1}) \times \dots \times \cP(M^{m_k})$.  Expand $\lambda$ to a $\tau$-structure by interpreting $F_\alpha^{n; m_1, \dots, m_k}$ as $F_\alpha$ restricted to the appropriate arity.  Then $D(M) \subset C(\bar{F})$.

Second, suppose we are given $M$ with universe $\lambda$.  WLOG, we can assume that $\tau$ has (second-order) Skolem functions.  Let $\bar{F} = \{F^M \mid F\in \tau\}$ be the functions of $M$ (interpreted as projection on other arities).  Then, since we have Skolem functions, $C(\bar{F}) \subset D(M)$.\hfill\dag

This model-theoretic characterization tends to be more useful to show the things that we want.

\begin{cor}
The superclub filter extends the club filter.
\end{cor}

{\bf Proof:}  The club filter can be characterized in the same way using first-order elementary substructure.  \hfill\dag\\

The next proposition shows that calling this a filter is justified.

\begin{prop}\label{sc-filter-prop}
The superclub filter is a $\kappa$-complete, possibly non-proper filter.
\end{prop}

{\bf Proof:}  The superclub filter is upwards closed by definition, so we only need to show it is closed under the intersection of $<\kappa$-many members.  Suppose that $X_\alpha \in \cF$ for $\alpha < \mu < \kappa$.  By Lemma \ref{sc-gen-lem}, there are $\tau_\alpha$-structures $M_\alpha$ with universe $\lambda$ and $|\tau_\alpha|<\kappa$ such that $D(M_\alpha) \subset X_\alpha$.  Set $\tau_*$ to be the disjoint union of the $\tau_\alpha$'s and $M_*$ to be the $\tau_*$-structure that is simultaneously an expansion of each $M_\alpha$.  Then $D(M_*) \subset D(M_\alpha)$, so this witnesses that $\bigcap_{\alpha<\mu} X_\alpha \in \cF$.

 It is a filter by definition.  Given $D(M_\alpha)$ with $M_\alpha$ a $\tau_\alpha$ structure, set $M_*$ to be the $\cup \tau_\alpha$-structure that is simultaneously an expansion of each $M_\alpha$.\hfill\dag\\

Proposition \ref{sc-filter-prop} leaves open the possibility that the superclub filter is non-proper, i.e., contains the empty set.  In fact, the properness of the superclub filter characterizes supercompact cardinals.

\begin{theorem}\label{superclub-thm}
Let $\kappa \leq \lambda$.
\begin{enumerate}
	\item If $\kappa$ {\bf is not} $\lambda$-supercompact, then there is an empty superclub.  Moreover, there is a uniform definition of the empty superclub.
	\item If $\kappa$ {\bf is} $\lambda$-supercompact, then every normal, fine, $\kappa$-complete ultrafilter on $\cP_\kappa \lambda$ extends the superclub filter.
\end{enumerate}
\end{theorem}

{\bf Proof:} If $\kappa$ is not $\lambda$-supercompact, then by Corollary \ref{2nd-lst-cor}, there is a structure $M$ of size $\lambda$ in a language $\tau$ of size $< \kappa$ such that $M$ has no $\prec_{\bL^2}$-substructures; looking at Magidor's proof, we can take $M$ to be some $(V_\beta, \in, \alpha)_{\alpha \in \mu}$ for some $\beta < \lambda$ and $\mu < \beta$.  WLOG $|M|$ is $\lambda$ by expanding the universe by a trivial sort.  Then $D(M) = \emptyset$ is in the superclub filter.

If $\kappa$ is $\lambda$-supercompact, then let $U$ be a normal, $\kappa$-complete ultrafilter on $\cP_\kappa \lambda$ and derive $j_U:V \to \cM_U$.  Let $M$ be a $\tau$-structure with universe $\lambda$ and, WLOG, it has Skolem functions.  By Magidor's result and its extensions, $j_U"M \subset_{\bL^2} j_U(M)$.  Since $M$ has Skolem functions, $j_U"M\prec_{\bL^2} j_U(M)$.  Thus, $j"\lambda \in \left(D(j_U(M))\right)^{\cM_U} = j\left(D(M)\right)$.  So $D(M) \in U$.  \hfill \dag\\

In particular, we get the following:

\begin{cor}\label{sc-sc-cor}
Given $\kappa \leq \lambda$, $\kappa$ is $\lambda$-supercompact iff the superclub filter on $\cP_\kappa \lambda$ is proper.
\end{cor}

While superclubs give a characterization of supercompact cardinals, they lack a nice characterization in the spirit of ``closed unbounded sets" that clubs have.  Such a characterization would shed light on properties of the $\prec_{\bL^2}$ relations and permit the exploration of superstationary sets.

\section{Extenders and omitting types} 
\label{extot-sec}

A key distinction between Theorem \ref{framework-thm} and Theorems \ref{second-thm} and \ref{sort-thm}  is the lack of an analogue of Theorem \ref{framework-thm}.(\ref{ult}) in the results of Section \ref{second-order-sec}.  The large cardinals of Section \ref{second-order-sec} are typically characterized by the existence of certain kinds of extenders, but the modifier ``certain" is typically characterized in a way that nakedly concerns embeddings between models of set theory--e.g., \cite[Exercise 26.7]{kanamori} characterizes $\kappa$ being $\lambda$-strong by the existence of an extender $E$ such that $j_E:V\to \cM_E$ witnesses $\lambda$-strength--rather than any combinatorial feature of the extender.

We begin with a combinatorial characterization, although it still references the $V_\alpha$'s and may be of limited interest.  However, we use this as a starting point to investigate a larger question: to what extent is second-order logic necessary to characterize the large cardinals in Section \ref{second-order-sec}?

We should say a few words about our definition of extenders.  We say that a \emph{$(\kappa, \lambda)$-extender} is $E = \{E_\ba \mid \ba \in [\lambda]^{<\omega}\}$, where each $E_\ba$ is a $\kappa$-complete ultrafilter on ${}^\ba \kappa$ satisfying coherence, normality, and well-foundedness (see \cite[Section 26]{kanamori} for a statement of these conditions in a slightly different formalism).  Note that this is a compromise between the original definition of Martin-Steel \cite{ms-extenders}--which took $\ba \in [V_\lambda]^{<\omega}$ and required $E_\ba$ to be on ${}^\ba V_\kappa$--and more modern presentations like \cite{kanamori}--which takes $\ba \in [\lambda]^{<\omega}$ but requires $E_\ba$ to be on $[\kappa]^{|\ba|}$.  Crucially,  we also depart from both of these and don't require that $\{s \in {}^\ba \kappa \mid \forall a_1, a_2 \in \ba \left(a_1 < a_2 \leftrightarrow s(a_1) < s(a_2)\right)\} \in E_\ba$.

Now we are ready to give a combinatorial characterization.

\begin{defin}
We say a bijection $h: \beth_\alpha \to V_\alpha$ is \emph{rank-layering} iff for every $\beta < \alpha$, $h"\beth_\beta = V_\beta$.
\end{defin}

Note this condition is equivalent to $x \in y \in V_\alpha$ implies $h^{-1}(x) \in h^{-1}(y)$.

\begin{prop} \label{strong-char-prop}
The following are equivalent.
\begin{enumerate}
	\item $\kappa$ is $\lambda$-strong.
	\item For \emph{some} rank-layering bijection $h:\beth_\lambda \to V_\lambda$, there is a $(\kappa, \beth_\lambda)$-extender $E$ such that for all $\alpha, \beta < \lambda$, we have
	$$h(\alpha) \in h(\beta) \iff \{ s \in {}^{\{\alpha, \beta\}}\kappa \mid h(s(\alpha)) \in h(s(\beta))\} \in E_{\{\alpha, \beta\}}$$
	\item For \emph{every} rank-layering bijection $h:\beth_\lambda \to V_\lambda$, there is a $(\kappa, \beth_\lambda)$-extender $E$ such that for all $\alpha, \beta < \lambda$, we have
	$$h(\alpha) \in h(\beta) \iff \{ s \in {}^{\{\alpha, \beta\}}\kappa \mid h(s(\alpha)) \in h(s(\beta))\} \in E_{\{\alpha, \beta\}}$$
\end{enumerate}
\end{prop}

In each of the cases, we have that $x \in V_\lambda$ is the image of $[\{h^{-1}(x)\}, s \mapsto h (s(h^{-1}(x)))]_E$ after the transitive collapse.

{\bf Proof of Proposition \ref{strong-char-prop}:} (3) implies (2) is clear.\\

$(2) \to (1)$: Given such an extender, set $j_E: V \to \cM$ to come from the extender power of $V$ by $E$ followed by the transitive collapse; then $\crit j = \kappa$ and $j(\kappa) > \lambda$ by standard arguments (see \cite[Lemmas 26.1 and 26.2]{kanamori}).  So we are left with showing that $V_\lambda \subset \cM$.  Let $\pi:\prod V/E \to \cM$ be the Mostowski collapse.  Following the above, we will show that, for each $x \in V_\lambda$, 
$$x = \pi \left([\{h^{-1}(x)\}, s \mapsto h (s( h^{-1}(x)))]_E\right)$$
We work by induction on the rank of $x$.  For each $y \in x$, the condition on our extender precisely gives that 
\begin{eqnarray*}
\{ s \in {}^{\{h^{-1}(y), h^{-1}(x)\}}\kappa& \mid& h(s(h^{-1}(y))) \in h(s(h^{-1}(x)))\} \in E_{\{h^{-1}(y), h^{-1}(x)\}}\\
\pi \left([\{h^{-1}(y)\}, s \mapsto h (s( h^{-1}(y)))]_E\right) &\in& \pi \left([\{h^{-1}(x)\}, s \mapsto h (s( h^{-1}(x)))]_E\right)
\end{eqnarray*}
Now suppose that $z \in \pi \left([\{h^{-1}(x)\}, s \mapsto h \circ s \circ h^{-1}(x)]_E\right)$.  We wish to show that $\pi^{-1}(z)$ is one of our terms for some $y \in x$.  We know that $\pi^{-1}(z)$ is $[\ba, f]_E$ for some $\ba \in [\beth_\lambda]^{<\omega}$ and $f$ with domain ${}^\ba\kappa$.  WLOG we may assume that $h^{-1}(x) \in \ba$, and coherence implies that $[\{h^{-1}(x)\}, s \mapsto h (s(h^{-1}(x)))]_E=[\ba, s \mapsto h (s(h^{-1}(x)))]_E$.  Thus, we have that 
$$\{s \in {}^\ba \kappa \mid f(s) \in h(s(h^{-1}(x)))\} \in E_\ba$$
If we set $\alpha = h^{-1}(x)$, then the rank-layering property tells us that 
$$\{s \in {}^\ba \kappa \mid h^{-1}(f(s)) \in s(\alpha)\} \in E_\ba$$
By normality, this means that there is $\bb \supset \ba$ and $\beta \in \bb$ such that
$$\{s \in {}^\bb \kappa \mid h^{-1}(f(s\rest \ba)) = s(\beta) \} \in E_\bb$$
If we set $y=h(\beta)$ and $g(s)=f(s\rest \ba)$, we can rewrite this as 
$$\{s \in {}^\bb \kappa \mid g(s) = h(s(h^{-1}(y)))\} \in E_\bb$$
Putting all of this together, we have that
$$\pi^{-1}(z) = [\ba, f]_E = [\bb, g]_E = [\{h^{-1}(y)\}, s \mapsto h(s(h^{-1}(y)))]_E$$
So $z$ must be an element of this form with $y \in x$, as desired.

$(1) \to (3)$: Fix $h:\beth_\lambda\to V_\lambda$ to be a rank-layering bijection and $j:V\to \cM$ to witness that $\kappa$ is $\lambda$-strong, so $\crit j = \kappa$, $j(\kappa)>\lambda$, and $V_\lambda \subset \cM$.  \cite[Chapter 26]{kanamori} describes the general method of deriving an extender $E^*$ from an elementary embedding by setting
$$X \in E^*_\ba \iff \id_\ba \in j(X)$$
The choice of the identity to `seed' the ultrafilters is not necessary; we can (and will) change this function without changing the fact that the derived ultrafilters form an extender.  In particular, define $E$ by specifying, for each $\ba \in [\beth_\lambda]^{<\omega}$ and $X \subset {}^\ba \kappa$, 
$$X \in E_\ba \iff j(h)^{-1} \circ h \circ j^{-1} \rest j(\ba) \in j(X)$$
Then $E$ is a $(\kappa, \beth_\lambda)$-extender and we only need to show the additional property.  We verify this with the following chain of equivalences, making crucial use of our choice of seed: for $\alpha, \beta < \lambda$, set $\ba = \{\alpha, \beta\}$
\begin{eqnarray*}
\{s \in {}^\ba\kappa \mid h(s(\alpha))\in h(s(\beta))\} &\in& E_\ba \iff j(h)^{-1} \circ h \circ j^{-1} \rest j(\ba) \in j\left( \{s \in {}^\ba\kappa \mid h(s(\alpha))\in h(s(\beta))\}\right)\\
&\iff& j(h)^{-1} \circ h \circ j^{-1} \rest j(\ba) \in \left\{s \in {}^{j(\ba)}j(\kappa) \mid j(h)\left(s(j(\alpha))\right)\in j(h)\left(s(j(\beta))\right)\right\}\\
&\iff& j(h)\left(j(h)^{-1}(h(j^{-1}(j(\alpha))))\right) =  j(h)\left(j(h)^{-1}(h(j^{-1}(j(\beta))))\right)\\
&\iff& h(\alpha) \in h(\beta)
\end{eqnarray*}
as desired. \hfill \dag\\

Now we turn to the issue of how necessary logics beyond $\bL_{\infty, \infty}$ are to characterize large cardinals and focus on strong cardinals.  Theorem \ref{strong-char-thm} characterizes strong cardinals in terms of the logic $\bL^2$; is there a characterization of strong cardinals solely in terms of $\bL_{\infty, \infty}$?

For $\kappa \leq \lambda$, consider the following (definable-class) $\bL_{\kappa, \omega}(Q^{WF})$-theory and types for $y \in V_\lambda$:
\begin{eqnarray*}
\tau &=& \{E, c_a, c, d_{a'}\}_{a \in V, a' \in V_\lambda}\\
T &=& ED_{\bL_{\kappa, \omega}(Q^{WF})}(V, \in, x)_{x \in V} \cup \{d_i = c_i < c < c_\kappa \mid i < \kappa\} \cup \{d_b \in d_a \mid b \in a \in V_\lambda\}\\
p_y(x) &=& \{ x E d_y \wedge \neg (x = d_z) \mid z \in y\}
\end{eqnarray*}
It follows from the methods of the previous sections that $T$ has a model omitting each $p_y$ iff $\kappa$ is $\lambda$-strong.  However, what we lack from Proposition \ref{elem-comp-thm} is an appropriate type omitting-compactness scheme and a decomposition of $T$ along that scheme that is satisfiable.  The difficulty in constructing this becomes more clear if we look at the extender product construction.

Given a $(\kappa,\lambda)$-extender $E$, for each $\ba \in [\lambda]^{<\omega}$, we form the ultraproduct $\prod V/E_{\ba}$.  The coherence axiom and \L o\'{s}' Theorem insures that the restriction function induces a coherent system of $\bL_{\kappa, \kappa}$-elementary embeddings  from $\prod V/E_\ba$ to $\prod V/E_\bb$ when $\ba \subset \bb$.  Since $[\lambda]^{<\omega}$ is directed, we can get a colimit $\prod V/E$ and $\bL_{\kappa, \omega}$-elementary embeddings.  Since $E$ is well-founded, so is $\prod V/E$, so we can form its transitive collapse $\cM_E$.

We can generalize this to a more general model-theoretic context.  Fix some $\ba \in [\lambda]^{<\omega}$ and suppose that we have a collection of $\tau$-structures $\{M_s \mid s \in {}^\ba \kappa\}$ and form the ultraproducts $\prod_{s \in {}^\bb \kappa} M_{s \rest \ba} /E_\ba$ for $\bb \supset \ba$.  Again, the coherence gives a coherent system of embeddings, so we can form the extender product $\prod M_s /E$ as the colimit of this system.  This structure has universe $\{[\bb, f]_E \mid \ba \subset \bb \in [\lambda]^{<\omega}, f \in {}^\bb \kappa\}$ where $[\bb, f]_E = [\bc, g]_E$ iff $\{s \in {}^{\bb\bc} \kappa \mid f(s \rest \bb) = g(s \rest \bc) \} \in E_{\bb\bc}$.  Then we have following result for \L o\'{s}' Theorem.

\begin{prop}\label{wf-los-prop}
Let $E = \{E_\ba \mid \ba \in [\lambda]^{<\omega}\}$ be a system of $\kappa$-complete ultrafilters satisfying coherence.  The following are equivalent:
\begin{enumerate}
	\item $E$ is well-founded.
	\item \L o\'{s}' Theorem holds for $\bL_{\kappa, \omega}(Q^{WF})$ formulas.  That is, given $\tau$-structures $\{M_s \mid s \in {}^\ba \kappa\}$, $\phi(x_1, \dots, x_n) \in \bL_{\kappa, \omega}(Q^{WF})(\tau)$, and $[\bb_1, f_1]_E, \dots, [\bb_n, f_n]_E \in \prod M_s/E$, we have
	\begin{center}
	$$\prod M_s/E \vDash \phi\left([\bb_1, f_1]_E, \dots, [\bb_n, f_n]_E\right)$$
iff
	$$\left\{ s \in {}^{\cup \bb_i} \kappa \mid M_{s \rest \ba} \vDash \phi\left( f_1(s\rest \bb_1), \dots, f_n(s \rest \bb_n)\right)\right\}\in E_{\cup \bb_i}$$
\end{center}
\end{enumerate}
\end{prop}

{\bf Proof:} For one direction, it is known that $E$ is well-founded iff $\prod V/E$ is well-founded, which follows from \L o\'{s}' Theorem applied to $Q^{WF}xy(x\in y)$.

For the other direction, fix $\bb \in [\lambda]^{<\omega}$ and $\tau$-structures $\{M_s \mid s \in {}^\bb \kappa\}$.  We show \L o\'{s}' Theorem for $\bL_{\kappa, \omega}(Q^{WF})$ by induction.  Standard arguments take care of everything but the $Q^{WF}$ quantifier.  So suppose \L o\'{s}' Theorem holds for $\phi(x, y, \bz)$ and $[\ba, f]_E \in M_E := \prod M_s/E$.

First, suppose that $\{s \in {}^\ba \kappa \mid M_{s\rest \bb} \vDash Q^{WF} xy \phi(x, y, f(s))\} \not\in E_\ba$.  Set $X$ to be the complement of this set.  For $s \in X$, there is $c^s_r \in M_{s \rest \bb}$ for $r < \omega$ such that $M_{s \rest \bb} \vDash \phi(c^s_{r+1}, c^2_r, f(s)$.  Then $c_r := [\ba, c^s_r]_E \in M_E$ witnesses the illfoundedness of $\phi$.

Second, suppose that $X_0 = \{s \in {}^\ba\kappa \mid M_{s \rest \bb} \vDash Q^{WF}xy\phi(x,y,f(s))\} \in E_\ba$ and $M_E \vDash \neg Q^{WF}xy\phi(x, y, [\ba, f]_E)$.  Then there is $[\ba_r, f_r]_E \in M_E$ such that
$$M_E \vDash \phi\left([\ba_{r+1}, f_{r+1}]_E, [\ba_r, f_r]_E, [\ba, f]_E\right)$$
for all $r < \omega$ and, WLOG, $\ba \subset \ba_r \subset \ba_{r+1}$.  Then
$$X_{r+1} := \{s \in {}^{\ba_{r+1}}\kappa \mid M_{s \rest \bb} \vDash \phi(f_{r+1}(s), f_r(s \rest \ba_r), f(s\rest \ba)\} \in E_{\ba_{r+1}}$$
The well-foundedness of $E$ gives $h:\cup \ba \to \kappa$ such that $h\rest \ba_{r+1} \in X_{r+1}$ and $h \rest \ba \in X_0$.  Then, for $r < \omega$,
$$M_{h\rest \bb} \vDash \phi(f_{r+1}(h\rest \ba_{r+1}), f_r(h \rest \ba_r), f(h\rest \ba))$$
Thus, $\seq{f_r(h \rest \ba_r) \in M_{h \rest \bb} \mid r < \omega}$ witnesses the illfoundedness of $\phi(x, y, f(h \rest \ba))$ in $M_{h \rest \bb}$, contradicting $h \rest \ba \in X_0$.\hfill \dag\\

So, to decompose $T$ above, one could try to construe $\prod V/ E$ as an extender product $\prod M_s /E$ in the appropriate language and for the appropriate $\ba$, and then see what fragment of $T$ $M_s$ satisfies.  However, the problem is that the factors that make up $\prod V/E$ don't have expansions to $\tau$-structures.  Rather, the analysis of $\prod V/E$ crucially uses that it can be seen as the extender power $\prod_{s\in {}^\ba\kappa} V/E$ for any choice of $\ba$.  Thus, there is no way to analyze which parts of the types each factor omits.

However, there is a nice criterion for when an extender product (or just a coherent ultraproduct by a $\kappa$-complete, well-founded coherent ultrafilter) omits a $\bL_{\kappa, \omega}(Q^{WF})$-type based on the behavior of the original models $M_s$ \emph{provided} that the domain of the type is a subset that appears as an element of $\cM_E$.

\begin{prop}
Let $E$ be a $(\kappa, \beth_\lambda)$-extender witnessing that $\kappa$ is $\lambda$-strong.  Suppose that $\ba \in [\beth_\lambda]^{<\omega}$, $\{M_s \mid s \in {}^\ba \kappa\}$ is a collection of $\tau$-structures, $A_\ell \subset \prod M_s/E$ of rank $\leq \lambda$ for $\ell = 0, 1$, and $\phi(x, y) \in \bL_{\kappa, \omega}(Q^{WF})$.  Set $p(x) = \{\phi(x, a) \mid a \in A_0\} \cup \{\neg \phi(x, a) \mid a \in A_1\}$ and $p_s(x) = \{\phi(x, a) \mid a \in f_0(s)\} \cup \{ \neg \phi(x, a) \mid a \in f_1(s)\}$, where $[\bb, f_\ell]_E$ represents $A_\ell$ in $\cM_E$.  Then, the following are equivalent:
\begin{enumerate}
	\item $\prod M_s/E$ omits $p$.
	\item $\{s \in {}^\bb \kappa \mid M_{s \rest \ba} \text{ omits } p_s\} \in E_{\bb}$.
\end{enumerate}
\end{prop}

Note that we have restricted both to the case of $\phi$-types and to the case where $E$ is an extender witnessing strength for simplicity and because those cases suffice for our assumption.  We could remove these assumptions, instead requiring that $A_\ell$ is an element of the subset sort of \mbox{$\prod (M_s, \cP(M_s), \in)/E$}.

{\bf Proof:} The structure $\prod M_s/E$ is (isomorphic to) $j_E(g)(\ba)$, where $g$ is the function taking $s \in {}^\ba \kappa$ to $M_s$ and $\ba = [\ba, s \mapsto s(\ba)]_E$.  Then, since $A_\ell \in V_\lambda \in \cM_E$, \L o\'{s}' Theorem for extenders tells us that 
\begin{eqnarray*}
\cM_E \vDash ``j_E(g)(\ba) \text{ omits } p" \iff \cM_E \vDash ``j_E(g)(\ba) \text{ omits the $\phi$-type generated by } A_0, A_1"\\
\iff \{s \in {}^\bb\kappa \mid g(s \rest \ba) \text{ omits the $\phi$-type generated by }f_0(s), f_1(s)\} \in E_\bb \\\iff \{s \in {}^\bb \kappa \mid M_{s \rest \ba} \text{ omits } p_s\} \in E_{\bb}
\end{eqnarray*}
\hfill \dag\\

Still, this doesn't give a syntactic characterization of strength because it deals with type omission for types over sets, whereas Theorem \ref{framework-thm} deals just with type omission over the empty set (in the appropriate language).  Thus, we are still left with the following question.

\begin{question}
Given $\kappa \leq \lambda$, is there a syntactic property of logics such that $\bL_{\kappa, \omega}(Q^{WF})$ (or some other sub-logic of $\bL_{\kappa, \kappa})$ satisfies this property iff $\kappa$ is $\lambda$-strong?
\end{question}

\bibliographystyle{amsalpha}
\bibliography{bib}
\end{document}